\documentclass[12pt,A4paper]{amsart}
\usepackage{amsmath,amsthm,amssymb,color}
\usepackage[bookmarksnumbered,colorlinks]{hyperref}
\usepackage{graphicx}
\hypersetup{colorlinks=true,linkcolor=red,anchorcolor=green,citecolor=cyan}

\newtheorem{theorem}{Theorem}[section]
\newtheorem{lemma}[theorem]{Lemma}
\newtheorem{proposition}[theorem]{Proposition}
\numberwithin{equation}{section}
\newcommand{\N}{\mathbb{N}}
\newcommand{\R}{\mathbb{R}}

\newcommand{\Q}{\mathbb{Q}}
\newcommand{\ve}{\varepsilon}

\begin{document}

\title[\tiny Almost uniform vs.~pointwise convergence from a linear point of view]{Almost uniform vs.~pointwise convergence from a linear point of view}

\author[Bernal]{L.~Bernal-Gonz\'alez}
\address[L. Bernal-Gonz\'alez]{\mbox{}\newline \indent Departamento de An\'alisis Matem\'atico \newline \indent Facultad de Matem\'aticas
\newline \indent Instituto de Matem\'aticas de la Universidad de Sevilla (IMUS)
\newline \indent Universidad de Sevilla
\newline \indent Avenida Reina Mercedes s/n, 41012-Sevilla (Spain).}
\email{lbernal@us.es}

\author[Calder\'on]{M.C.~Calder\'on-Moreno}
\address[M.C.~Calder\'on-Moreno]{\mbox{}\newline \indent Departamento de An\'alisis Matem\'atico \newline \indent Facultad de Matem\'aticas
\newline \indent Instituto de Matem\'aticas de la Universidad de Sevilla (IMUS)
\newline \indent Universidad de Sevilla
\newline \indent Avenida Reina Mercedes s/n, 41012-Sevilla (Spain).}
\email{mccm@us.es}

\author[Gerlach]{P.J.~Gerlach-Mena}
\address[P.J.~Gerlach-Mena]{\mbox{}\newline \indent Departamento de Estad\'{\i}stica e Investigaci\'on Operativa
\newline \indent Facultad de Matem\'aticas
\newline \indent Universidad de Sevilla
\newline \indent Avenida Reina Mercedes s/n, 41012-Sevilla (Spain).}
\email{gerlach@us.es}

\author[Prado]{J.A.~Prado-Bassas}
\address[J.A.~Prado-Bassas]{\mbox{}\newline \indent Departamento de An\'alisis Matem\'atico
\newline \indent Facultad de Matem\'aticas
\newline \indent Instituto de Matem\'aticas de la Universidad de Sevilla (IMUS)
\newline \indent Universidad de Sevilla
\newline \indent Avenida Reina Mercedes s/n, 41012-Sevilla (Spain).}
\email{bassas@us.es}



\begin{abstract}
A review of the state of the art of the comparison between any two different modes of convergence of sequences of measurable functions is carried out with focus on the algebraic structure of the families under analysis. As a complement of the amount of results obtained by several authors, it is proved, among other assertions and under natural assumptions, the existence of large vector subspaces as well as of large algebras contained in the family of the sequences of measurable functions converging to zero pointwise almost everywhere but not almost uniformly, and in the family of the sequences of measurable functions converging to zero almost uniformly but not uniformly almost everywhere.
\end{abstract}

\maketitle

\section{Introduction}

\quad Several kinds of convergence have been considered in the context of the class of sequences of measurable functions defined on a positive measure space. It is well known that, except for very few cases, these modes of convergence are pairwise different. This means that, if \,$A, \, B$ \,are two distinct modes of convergence such that \,$A$ \,implies \,$B$, then one can find a sequence of measurable functions that is $B$-convergent but not $A$-convergent.
The relationships between these diverse modes of convergence is studied in many books, as for instance \cite{billingsley,bogachev,fremlin,karr,nielsen,stoyanov,wisehall}, among others.
For background on measure theory and convergence, we refer the interested reader to them. In this paper it will be analyzed to what extent --in a linear sense to be specified later-- the failure of the mentioned reverse implications takes place.

\vskip 3pt

Let us set the basic notation and recall the main concepts. Assume that \,$(\Omega,\mathcal{A},\mu )$ \,is a positive measure space, and denote by \,$L_0$ \,the
vector space of all ($\mu$-classes of) measurable functions \,$f : \Omega \to \R$. Hence, two functions \,$f,g \in L_0$ \,are equal if
they are equal almost everywhere, that is, if
there is a measurable set \,$Z \in \mathcal{A}$ \,with \,$\mu (Z) = 0$ \,such that \,$f(x) = g(x)$ \,for all \,$x \in \Omega \setminus Z$. If \,$\mu (\Omega ) = 1$ \,then
\,$(\Omega,\mathcal{A},\mu )$ \,is usually called a probability space, and a measurable function on it is called a random variable. Let us consider the
vector space \,$L_0^{\N}$ \,of all sequences of measurable functions on \,$\Omega$. Assume that \,$f \in L_0$ \,and \,$(f_n) \in L_0^{\N}$. Then we say that
\begin{enumerate}
\item[$\bullet$] $(f_n)$ \,tends to \,$f$ \,{\it pointwise almost everywhere} provided that there exists a set \,$Z \in \mathcal{A}$ \,with \,$\mu (Z) = 0$
 \,such that \,$f_n(x) \to f(x)$ \,for every \,$x \in \Omega \setminus Z$.
\item[$\bullet$] $(f_n)$ \,tends to \,$f$ \,{\it uniformly almost everywhere} provided that there exists a set \,$Z \in \mathcal{A}$ \,with \,$\mu (Z) = 0$
 \,such that \,$f_n \to f$ \,uniformly on \,$\Omega \setminus Z$.
\item[$\bullet$] $(f_n)$ \,tends to \,$f$ \,{\it in measure} whenever \,$\displaystyle{\lim_{n \to \infty}} \mu (\{ x \in \Omega : \, |f_n(x) - f(x)| > \varepsilon \}) = 0$
for every \,$\varepsilon > 0$.
\item[$\bullet$] $(f_n)$ \,tends to \,$f$ \,{\it almost uniformly} if, for every \,$\varepsilon > 0$, there is a set \,$Z_\varepsilon \in \mathcal{A}$ \,with
                \,$\mu (Z_\varepsilon) < \varepsilon$ \,such that \,$f_n \to f$ \,uniformly on \,$\Omega \setminus Z_\varepsilon$.
\item[$\bullet$] $(f_n)$ \,tends to \,$f$ \,{\it in $q$-norm} (or {\it in $q$-mean}), where \,$q \in (0,+\infty )$, if \break
$\displaystyle{\lim_{n \to \infty}} \int_{\Omega} |f_n - f|^q \,d\mu = 0$.
\item[$\bullet$]$(f_n)$ \,tends to \,$f$ \,{\it completely} \,provided that, for any prescribed \,$\varepsilon > 0$,
the series \,$\sum_{n=1}^{\infty} \mu (\{ x \in \Omega : \, |f_n(x) - f(x)| > \varepsilon \})$ \,converges.
\end{enumerate}

Note that the six given definitions are consistent under replacement of \,$f$ \,and each \,$f_n$ \,by respective a.e.-equal measurable functions \,$\widetilde{f}$, $\widetilde{f_n}$ (a.e.~stands for ``almost everywhere''). There are a number of implications among these kinds of convergence. In all of them, the limit function \,$f$ \,in the conclusion is the same as the one in the assumption. We next list them, taking into account that some of the properties happen only in finite measure spaces. The proofs in this case are usually given in probability spaces, to which case it is always possible to reduce the problem just by dividing into \,$\mu (\Omega )$. The reader is referred to the books cited at the beginning, and especially to \cite[Chap.~22]{nielsen}, to check these implications:
\begin{enumerate}
\item[$\bullet$] Uniform convergence a.e. \,$\Longrightarrow$ \, almost uniform convergence.
\item[$\bullet$] Almost uniform convergence \, $\Longrightarrow$ \, pointwise convergence a.e.~and convergence in measure.
\item[$\bullet$] For \,$q \in (0,+\infty )$, convergence in $q$-mean \, $\Longrightarrow$ \, convergence in measure.
\item[$\bullet$] Complete convergence \, $\Longrightarrow$ \, convergence in measure and pointwise convergence a.e.
\item[$\bullet$] If \,$\mu$ \,is finite, then uniform convergence a.e. \,$\Longrightarrow$ \, complete convergence and convergence in $q$-mean for every $q > 0$.
\item[$\bullet$] (Egoroff's theorem) If \,$\mu$ \,is finite, then pointwise convergence a.e. $\Longrightarrow$ \, almost uniform convergence.
                 Hence, both modes of convergence are equivalent in this context.
\end{enumerate}

As complementary results, we have that convergence in measure entails almost uniform convergence of a subsequence.
Moreover, if $q \in [1,+\infty)$, then pointwise convergence a.e.~plus existence of \,$g \in L^q$ \,with \,$|f_n| \le g$ a.e.~for all $n \in \N$ implies convergence in $q$-mean
(Lebesgue dominated convergence theorem).

\vskip 3pt

Since \,$f_n \to f$ \,(in any of the given modes) if and only if \,$f_n - f \to 0$, the study of the comparison of kinds of convergence can be reduced to study how many sequences
\,$(f_n) \in L_0^{\N}$ \,satisfy \,$f_n \to 0$ \,in one sense but not in other sense. In this framework, several authors in a number of papers \cite{araujobmps,vecina,polacos2018,polacos2023,bernalordonez2014,calderongerlachpradomodes,fesetu} have invested much effort in the analysis of the algebraic size or even the algebraic-topological size of the families of these special sequences.

\vskip 3pt

Our aim in this paper is to contribute to go deeper into this analysis, with focus on the class of the sequences of measurable functions converging
to zero pointwise almost everywhere but not almost uniformly, and on the class of the sequences of measurable functions converging to zero almost
uniformly but not uniformly almost everywhere. Comparison with convergence to zero in measure, in $q$-mean, and complete convergence, will be also studied.
This will be carried out and precisely formulated in sections 3 and 4, while section 2
will be devoted to fix some pertinent notation, as well as to state the necessary algebraic and topological background and to recall the state
of the art of the subject of the paper.

\section{Notation and preliminaries}

\quad Firstly, we proceed by fixing some terminology coming from the modern theory of lineability.
This is a subject in functional analysis that has been thoroughly investigated for the last three decades.
For an extensive study of this theory, the reader is referred to the book \cite{ABPS} (see also \cite{bams,seoanetesis}).
The goal of lineability is to find large linear structures inside nonlinear subsets of a vector space.

\vskip 3pt

A subset \,$A$ \,of a vector space \,$X$ \,is called {\it lineable} whenever there is an infinite dimensional vector subspace of \,$X$ \,that is contained, except for zero, in \,$A$; and \,$A$ \,is said to be {\it algebrable} if it is contained in some linear algebra and there is an infinitely generated algebra contained, except for zero, in \,$A$.  More precisely, if \,$\alpha$ \,is a cardinal number, then \,$A$ \,is said to be {\it $\alpha$-lineable} if there exists a vector subspace \,$V \subset X$ \,such that \,$\dim(V) = \alpha$ \,and \,$V \setminus \{0\} \subset A$. In addition, if \,$A$ \,is contained in some commutative linear algebra, then \,$A$ \,is called {\it strongly $\alpha$-algebrable} if there exists an $\alpha$-generated free algebra \,$M$ \,with \,$M \setminus \{0\} \subset A$; that is, there exists a subset \,$B$ \,of cardinality \,$\alpha$ \,with the following property: for any positive integer \,$N$, any nonzero polynomial \,$P$ \,in \,$N$ \,variables without constant term and any distinct elements \,$x_1, \ldots ,x_N \in B$, we have
\begin{equation*}
  P(x_1, \ldots ,x_N) \in A \setminus \{0\}.
\end{equation*}
Plainly, strong $\alpha$-algebrability implies $\alpha$-lineability.
Now, assume that \,$X$ \,is a to\-po\-lo\-gi\-cal vector space and \,$A \subset X$. Then \,$A$ \,is called {\it dense-lineable} in \,$X$ \,if there is a dense vector subspace
\,$V \subset X$ \,such that \,$V \setminus\{0\} \subset A$; if, in addition, $\dim(V) = \alpha$, then \,$A$ \,is said to be {\it $\alpha$-dense-lineable}. Finally, if there is a closed infinite dimensional subspace \,$V \subset X$ \,such that \,$V \setminus \{0\} \subset A$, then \,$A$ \,is called {\it spaceable}.

\vskip 3pt

In order to apply the preceding concepts to our sequence space $L_0^{\N}$, a natural topology is needed in it.
With this aim, recall that a sequence \,$(f_n) \in L_0^{\N}$ \,is said {\it to converge locally in measure} to \,$f$ \,if,
for every \,$\ve > 0$ \,and every \,$A \in \mathcal A$ \,with \,$\mu (A) < \infty$, one has
$$
\lim_{n \to \infty} \mu (\{x \in A : \, |f_n(x) - f(x)| > \ve   \}) = 0.
$$
Then the family of sets
$$
\Big\{ f \in L_0 : \, \int_K \min \{1,|f|\} \, d\mu < \ve \Big\} \ \  \big( \ve > 0; \ K \in \mathcal{A} \hbox{ with }  \mu (K) < \infty \big)
$$
forms a base of $0$-neighborhoods for a topology on \,$L_0$ \,which makes it a topological vector space, so that convergence in this topology is
local convergence in measure. We have that this topological vector space \,$L_0$ \,is metrizable if and only if \,$(\Omega , \mathcal{A}, \mu)$ \,is $\sigma$-finite, that is, if and only if there is a sequence \,$(K_n) \subset \mathcal{A}$ \,such that \,$\mu (K_n) < +\infty$ \,for all \,$n \in \N$ \,and \,$\Omega = \bigcup_{n \in \N} K_n$
(see \cite[Chapter 24]{fremlin}). Of course, convergence in measure and local convergence in measure coincide if \,$\mu$ \,is finite.
Along this paper, the topology that we shall consider on \,$L_0$ \,is the one of {\it local convergence in measure,} while \,$L_0^{\N}$ \,will be endowed with
the corresponding {\it product topology.}

\vskip 4pt

Taking into account the object of our study, it is at this point convenient to introduce the following notation for a number
of families of sequences of measurable functions:

\vskip 4pt

$\mathcal{S}_{p} = \{ {\bf f} = (f_n) \in L_0^{\N}: \, f_n \to 0 \text{ pointwise almost everywhere} \}$

\vskip 3.5pt

$\mathcal{S}_{u} = \{ {\bf f} = (f_n) \in L_0^{\N}: \, f_n \to 0 \text{ \,uniformly almost everywhere} \}$

\vskip 3.5pt

$\mathcal{S}_{au} = \{ {\bf f} = (f_n) \in L_0^{\N}: \, f_n \to 0 \text{ \,almost uniformly} \}$

\vskip 3.5pt

$\mathcal{S}_{m} = \{ {\bf f} = (f_n) \in L_0^{\N}: \, f_n \to 0 \text{ \,in measure} \}$

\vskip 3.5pt

$\mathcal{S}_{L_q} = \{ {\bf f} = (f_n) \in L_0^{\N}: \, f_n \to 0 \text{ \,in $q$-mean} \} \ \ \ (q \in (0,+\infty ))$

\vskip 3.5pt

$\mathcal{S}_{c} = \{ {\bf f} = (f_n) \in L_0^{\N}: \, f_n \to 0 \text{ \,completely} \}$.

\vskip 5pt

Next, we try to gather the results obtained up to date concerning the failure of implications between the diverse modes of convergence {\it to zero.}
By \,$\mathfrak{c}$ \,we shall denote, as usual, the cardinality of the continuum, and by \,$\lambda$ \,the Lebesgue measure on
a Borel subset of \,$\R$. To avoid trivialities, it will always be assumed that \,$\mu (\Omega ) > 0$. Recall that if \,$\Omega$ \,is a topological space then
its Borel family \,$\mathcal{B}$ \,is the $\sigma$-algebra generated by the set of open subsets of \,$\Omega$.
A measure space \,$(\Omega ,\mathcal{A},\mu )$ \,is said to be {\it nonatomic} if it lacks atoms; recall that
a measurable set \,$A \in \mathcal{A}$ \,is said to be an atom if there do {\it not} exist \,$B,C \in \mathcal{A}$ \,such that
$$
B \cap C = \varnothing, \, B \cup C = A  \quad \hbox{and}  \quad \mu (B) > 0 < \mu (C).
$$
Under the previous notation and terminology, the following
facts have been established:
\begin{enumerate}
\item[(A)] In \cite[Theorem 7.1]{araujobmps} Ara\'ujo {\it et al.}~proved that \,$\mathcal{S}_m \setminus \mathcal{S}_p$ \,is $\mathfrak{c}$-lineable in the case
\,$([0,1],\mathcal{B},\lambda )$. This was improved by Calder\'on {\it et al.}~in \cite[Theorems 2.2 and 2.3]{calderongerlachpradomodes}, where it has been
respectively established its strong $\mathfrak{c}$-algebrability and its spaceability in $L_0^{\N}$. It is worth mentioned that the set considered in
\cite{calderongerlachpradomodes} was \,$\mathcal{S}_m \setminus (\mathcal{S}_p \cup \mathcal{S}_{au})$, but this set equals \,$\mathcal{S}_m \setminus \mathcal{S}_p$
\,(or \,$\mathcal{S}_m \setminus \mathcal{S}_{au}$)
\,thanks to Egoroff's theorem. Recently, Bartoszewicz {\it et al.}~in \cite[Theorem 11]{polacos2023} (see also their Theorem 13) succeeded in strengthening the algebrability part for any nonatomic probability space. In fact, together with other interesting additional properties, it can be derived from their results that even
\,$\mathcal{S}_m \setminus (\mathcal{S}_p \cup \mathcal{S}_{L_1})$ \,is strongly $\mathfrak{c}$-algebrable in this case.
\item[(B)] In \cite[Theorem 3.3]{calderongerlachpradomodes} it is shown that, in the case \,$([0,1],\mathcal{B},\lambda )$, the set
\,$(\mathcal{S}_p \cap \mathcal{S}_{au}) \setminus \mathcal{S}_u$ \,(again, Egoroff's theorem indicates that this set equals
\,$\mathcal{S}_p \setminus \mathcal{S}_u$ \, or \,$\mathcal{S}_{au} \setminus \mathcal{S}_u$) is strongly $\mathfrak{c}$-algebrable, as well as spaceable and $\mathfrak{c}$-dense lineable in \,$L_0^{\N}$.
\item[(C)] The first conclusion in the preceding item was improved by Ara\'ujo {\it et al.}, who in \cite[Theorem 1]{vecina} proved that the smaller set
\,$\mathcal{S}_{c} \setminus \mathcal{S}_u$ \,is strongly $\mathfrak{c}$-algebrable in the case \,$([0,1],\mathcal{B},\lambda )$.
In a general nonatomic probability space \,$(\Omega , \mathcal{A},\mu )$, they were able to show \cite[Theorem 8]{vecina} that the latter set is $\mathfrak{c}$-lineable.
\item[(D)] In \cite[Theorem 9]{vecina} it is demonstrated that \,$\mathcal{S}_p \setminus \mathcal{S}_c$ \,is $\mathfrak{c}$-lineable in a nonatomic probability space
\,$(\Omega, \mathcal{A},\mu )$. If, in addition, $\Omega$ \,is a complete separable metric space and \,$\mathcal{A}$ \,is the Borel $\sigma$-algebra in \,$\Omega$, then \,$\mathcal{S}_p \setminus \mathcal{S}_c$ \,is even strongly $\mathfrak{c}$-algebrable \cite[Theorem 2 and Remark 3]{vecina}.
\item[(E)] In Theorems 10, 11 and 12 of \cite{vecina} it is proved, respectively, the $\mathfrak{c}$-lineability of
\,$\left(\bigcap_{q \in (0,p)} \mathcal{S}_{L_q}  \right) \setminus \mathcal{S}_{L_p}$ (for each $p > 0$),
$\mathcal{S}_p \setminus \left( \bigcup_{q > 0} \mathcal{S}_{L_q} \right)$ \,and \,$\left( \bigcap_{q > 0} \mathcal{S}_{L_q} \right) \setminus \mathcal{S}_p$,
where \,$(\Omega , \mathcal{A},\mu )$ \,is assumed to be a nonatomic probability space. Four years before,
Fern\'andez-S\'anchez {\it et al.}~had established a result \cite[Theorem 2.4]{fesetu} which yielded as a consequence the strong $\mathfrak{c}$-algebrability of
\,$\mathcal{S}_p \setminus \mathcal{S}_{L_1}$, provided that the probability space \,$(\Omega , \mathcal{A},\mu )$ \,is not isomorphic
to any probability space \,$(\{1,2, \dots ,N\},2^{\{1,2, \dots ,N\}},\nu )$ ($N \in \N$). Note that, since \,$\mathcal{S}_{L_1} \subset \mathcal{S}_m$, it follows
from \cite[Theorem 12]{vecina} that \,$\mathcal{S}_m \setminus \mathcal{S}_p$ \,is strongly $\mathfrak{c}$-algebrable in the case of a nonatomic probability space,
so improving \cite[Theorem 2.2]{calderongerlachpradomodes}
mentioned in (A) above.
\item[(F)] In \cite{calderongerlachpradomodes} it is established, for the measure space $([0,+\infty ), \mathcal{B},\lambda )$, that the set
\,$\mathcal{S}_u \setminus \mathcal{S}_{L_1}$ is strongly $\mathfrak{c}$-algebrable (Theorem 4.1), and that its subset
\,$(\mathcal{S}_u \cap L_1^{\N}) \setminus \mathcal{S}_{L_1}$ \,is $\mathfrak{c}$-dense-lineable and spaceable in \,$L_1^{\N}$ (Theorem 4.2).
Here the Lebesgue space \,$L_1 = L_1([0,+\infty ))$ \,carries the usual $1$-norm topology, while \,$(L_1([0,+\infty )))^{\N}$
\,is endowed with the associated product topology.
\end{enumerate}
We emphasize the fact that the results stated for probability spaces hold for finite measure spaces.

\vskip 3pt

Sometimes it will be enough the assumption of a semifinite measure, which is a strictly weaker notion than $\sigma$-finiteness (see \cite[p.~141]{nielsen}).
Recall that if  \,$(\Omega , \mathcal{A},\mu )$ \,is a measure space, then both it and the measure \,$\mu$ \,are said to be
{\it semifinite} if \,$\mu (A) = \sup \{ \mu (B): \, B \in \mathcal{A}, \, B \subset A, \hbox{ and } \mu (B) < \infty \} \, \hbox{ for each } \,A \in \mathcal{A}$.
For instance, if \,$\Omega$ \,is an uncountable set then the counting measure on the $\sigma$-algebra \,$\mathcal{P}(\Omega )$ \,of all parts of \,$\Omega$ \,is semifinite but not $\sigma$-finite.
The next useful result can be found in \cite[Propositions 11.12 and 11.27]{nielsen}.

\begin{proposition} \label{Prop-caracterization semifinite, and atomless}
The measure space \,$(\Omega , \mathcal{A},\mu )$ \,is semifinite if, and only if, for each set \,$A \in \mathcal{A}$ \,with \,$\mu (A) > 0$
\,there is a set \,$B \in \mathcal{A}$ \,such that \,$B \subset A$ \,and \,$0 < \mu (B) < \infty$. If, in addition, the measure space is nonatomic,
then for every \,$A \in \mathcal{A}$ \,and every \,$\sigma \in [0,\mu (A)]$ \,there is \,$S \in \mathcal{A}$ \,such that \,$S \subset A$ \,and \,$\mu (S) = \sigma$.
\end{proposition}

The following two assertions will be invoked to find dense lineability of several families of sequences of measurable functions.
The first of them is a rather general criterion and can be extracted from \cite[Section 7.3]{ABPS}, while the second one is
more specific and is due to Ara\'ujo {\it et al.}~\cite[Theorem 6]{vecina}; it will be stated in a more set-theoretical way than in \cite{vecina}.

\begin{theorem}\label{A stronger than B}
Let \,$\alpha$ \,be an infinite cardinal number. Assume that \,$X$ is a metrizable separable topological vector space, and
that \,$A$ \,and \,$B$ \,are subsets of \,$X$ fulfilling the following properties:
\begin{enumerate}
  \item[\rm (i)] $A + B\subset A$.
  \item[\rm (ii)] $A \cap B = \varnothing$.
  \item[\rm (iii)] $A$ is $\alpha$-lineable.
  \item[\rm (iv)] $B$ is dense-lineable.
\end{enumerate}
Then \,$A$ \,is $\alpha$-dense-lineable.
\end{theorem}

When applying Theorem \ref{A stronger than B}, we shall frequently make use of the fact that the set
$$
c_{00}(L_0) := \{ {\bf f} = (f_n) \in L_0^{\N} : \hbox{ there exists } k = k({\bf f}) \in \N \hbox{ such that } f_n = 0 \hbox{ for all } n > k \}
$$
is a dense vector subspace of \,$L_0^{\N}$, hence dense-lineable in $L_0^{\N}$.

\begin{theorem}\label{Tma de Vecina dense-lineable}
Let \,$V$ \,be a vector space over \,$\R$, and \,$\mathcal{C}, \mathcal{D}$ \,two subsets of \,$V^{\N}$ \,satisfying the following properties, where \,$\mathcal{E}$
\,denotes any of such two subsets:
\begin{enumerate}
\item[\rm (a)] $\alpha \,\mathcal{E} \subset \mathcal{E}$ \,for all \,$\alpha \in \R$.
\item[\rm (b)] The zero sequence \,${\bf 0} := (0,0,0, \dots ) \in \mathcal{E}$.
\item[\rm (c)] If a sequence \,${\bf f} \in V^{\N}$ \,can be decomposed into finitely many subsequences belonging to \,$\mathcal{E}$,
then \,${\bf f} \in \mathcal{E}$.
\item[\rm (d)] If \,${\bf f} \in \mathcal{E}$, then all its subsequences also belong to \,$\mathcal{E}$.
\item[\rm (e)] If \,${\bf f} \in \mathcal{E}$ \,and finitely many terms are deleted from {\rm (}or added to{\rm )} \,${\bf f}$,
               then the new sequence also belongs to \,$\mathcal{E}$.
\item[\rm (f)] $\mathcal{C} \setminus \mathcal{D} \ne \varnothing$.
\end{enumerate}
Then \,$\mathcal{C} \setminus \mathcal{D}$ \,is $\mathfrak{c}$-lineable.
\end{theorem}

Note that each of the families \,$\mathcal{S}_p, \,\mathcal{S}_m, \mathcal{S}_u, \,\mathcal{S}_c, \,\mathcal{S}_{au}, \,\mathcal{S}_{L_q}$
\,satisfies all properties (a) to (e) of Theorem \ref{Tma de Vecina dense-lineable}.


\section{Dense lineability and spaceability}

\quad In this section and the next one we extend some results taken from those mentioned in section 2, and provide some new ones.

\vskip 3pt

Recall that we are assuming along this paper that \,$L_0$ \,carries the topology of local convergence in measure, while \,$L_0^{\N}$ \,carries the associated product topology. In order to study the denseness of a number of subsets of \,$L_0^{\N}$, it is convenient to have at our disposal a natural condition guaranteeing the separability of this space. With this aim, we consider the following condition that may or may not be satisfied by a measure space \,$(\Omega, \mathcal{A}, \mu )$:

\vskip 3pt

\noindent \phantom{aaa} (S) {\it There is a countable family \,$\mathcal A_1 \subset \mathcal A$ \,such that, for given \,$M \in \mathcal A$ \,with \hfil\break
\phantom{aaaaaa\,} $\mu (M) < \infty$ \,and \,$\ve > 0$, there exists \,$M_0 \in \mathcal A_1$ \,with \,$\mu (M \Delta M_0) < \ve$.}

\vskip 3pt

As usual, we have denoted $A \Delta B := (A \setminus B) \cup (B \setminus A)$, the symmetric difference of two sets $A,B$.
Examples of measure spaces satisfying (S) are the counting measure on \,$\N$ \,(take \,$\mathcal A_1 = \{$finite parts of \,$\N \}$) and the Lebesgue measure on
\,$\R^N$ \,(select as \,$\mathcal A_1$ \,the collection of all finite unions of open rectangles whose vertices have rational coordinates).

\begin{lemma} \label{Lemma L0N metrizable and separable}
Assume that \,$(\Omega, \mathcal{A}, \mu )$ is a $\sigma$-finite measure space satisfying condition {\rm (S).} Then \,$L_0^{\N}$ \,is metrizable and separable.
\end{lemma}

\begin{proof}
It is well known (see, e.g., \cite{engelking}) that a countable product of metrizable (se\-pa\-ra\-ble, resp.) topological spaces is metrizable (separable, resp.).
Hence, it suffices to show that \,$L_0$ \,is metrizable and separable. Its metrizability follows from the assumption of $\sigma$-finiteness of the measure space (see section 1).
It remains to prove the separability of \,$L_0$. Since convergence in $1$-mean implies convergence in measure, hence local convergence in measure, it is enough to check (i) the separability of the subspace
\,$L_1 = \{ f \in L_0 : \, \|f\|_1 < +\infty \}$ \,under the norm \,$\|f\|_1 = \int_\Omega |f| \, d\mu$ \,and (ii) the denseness of \,$L_1$ \,in \,$L_0$.\
Finally, (i) and (ii) are well known results: for (i) we refer to, for instance, \cite{ricker} or \cite[p.~137]{zaanen} (where it is proved that (S) is, in fact, a necessary and sufficient condition for the separability of $L_1$), while a detailed proof of (ii) can be found in \cite[Lemma 2]{antifubini} (condition (S) is not needed here).
\end{proof}

\begin{figure}[h!]\label{figura}
  \includegraphics[width=.22\textwidth]{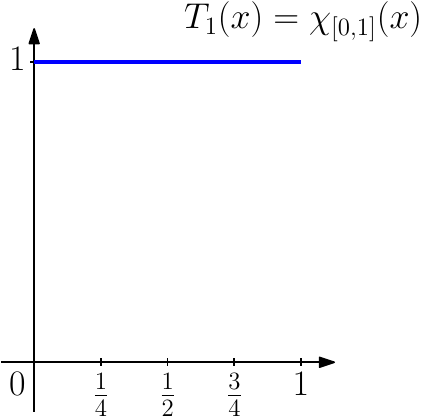}\\
  \includegraphics[width=.22\textwidth]{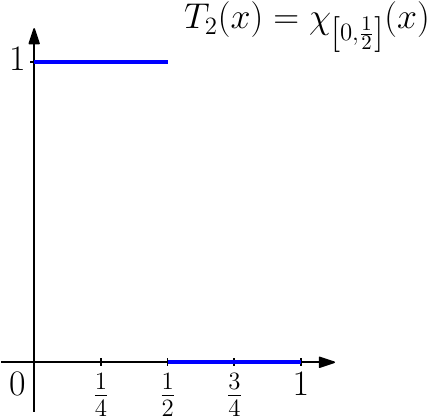}\quad
  \includegraphics[width=.22\textwidth]{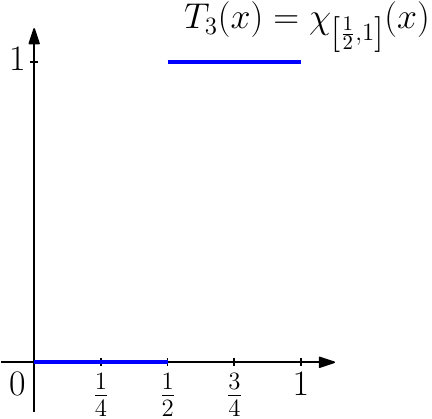}\\
  \includegraphics[width=.22\textwidth]{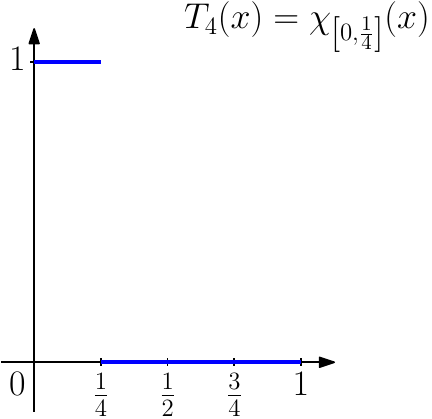}\quad
  \includegraphics[width=.22\textwidth]{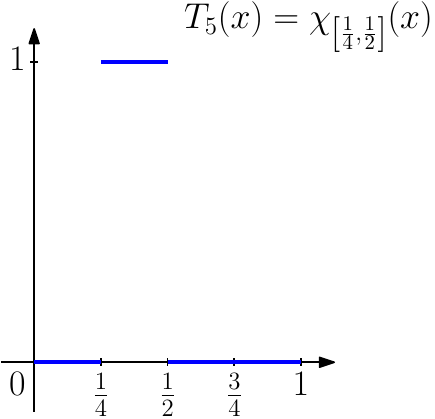}\quad
  \includegraphics[width=.22\textwidth]{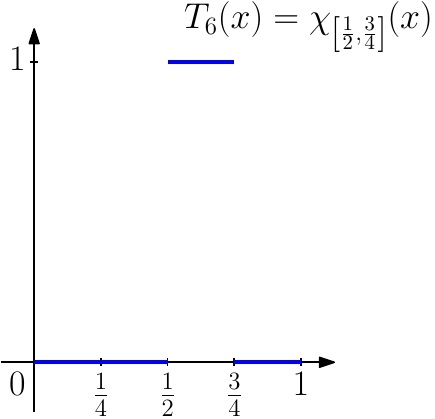}\quad
  \includegraphics[width=.22\textwidth]{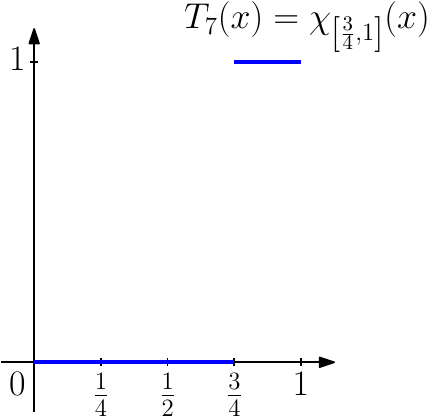}
  \caption{First iterations of the typewriter sequence.}
\end{figure}

\begin{theorem} \label{Sm-Spae dense lineable}
Let \,$(\Omega ,\mathcal{A},\mu )$ \,be a nonatomic $\sigma$-finite measure space satisfying condition {\em (S)}.
Then \,$\displaystyle{ \Big( \bigcap_{q > 0} \mathcal{S}_{L_q} \Big)} \setminus \mathcal{S}_{p}$ \,is $\mathfrak{c}$-dense-lineable in \,$L_0^{\N}$.
In particular, the sets \,$\mathcal{S}_m \setminus \mathcal{S}_{p}$, $\mathcal{S}_m \setminus \mathcal{S}_{c}$, $\mathcal{S}_m \setminus \mathcal{S}_{au}$ \,and \,$\mathcal{S}_m \setminus \mathcal{S}_{u}$
\,are $\mathfrak{c}$-dense-lineable in \,$L_0^{\N}$.
\end{theorem}

\begin{proof}
Let us first prove that  \,$\displaystyle{ \Big( \bigcap_{q > 0} \mathcal{S}_{L_q} \Big)} \setminus \mathcal{S}_{p}$ \,is $\mathfrak{c}$-lineable.
This is a mere adaptation of the proof of Theorem 12 in \cite{vecina}, but we shall do it for the sake of completeness.
In fact, the construction is an extension of the so-called ``typewriter sequence'' $(T_n)$ in the interval $[0,1]$ (see Figure 1). 
Since $\sigma$-finiteness implies semifiniteness, Proposition \ref{Prop-caracterization semifinite, and atomless} yields the existence of sets
\,$A,S \in \mathcal{A}$ \,such that \,$0 < \mu (A) < \infty$, $S \subset A$, and \,$\mu(S) = {\mu (A) \over 2}$.
Defining \,$\Omega_0 := S$ \,and \,$\Omega_1 := A \setminus S$, we get
\,$\Omega_0 \cup \Omega_1 = A$, $\Omega_0 \cap \Omega_1 = \varnothing$, and \,$\mu (\Omega_0) = {\mu (A) \over 2} = \mu (\Omega_1)$.
Applying the same argument to each of the sets \,$\Omega_0, \, \Omega_1$, and repeating it recursively, we obtain a family
\,$\big\{ \Omega_{j_1 \cdots j_m} : \, m \in \N, \, j_1, \dots ,j_m \in \{0,1\} \big\} \subset \mathcal{A}$ \,satisfying the following properties
for each \,$m \in \N$:
\begin{enumerate}
\item[$\bullet$] $\Omega_{j_1 \cdots j_m} = \Omega_{j_1 \cdots j_m 0} \cup \Omega_{j_1 \cdots j_m 1}$.
\item[$\bullet$] $\Omega_{j_1 \cdots j_m 0} \cap \Omega_{j_1 \cdots j_m 1} = \varnothing$.
\item[$\bullet$] $\mu (\Omega_{j_1 \cdots j_m 0}) =  \mu(\Omega_{j_1 \cdots j_m 1}) = {\mu (A) \over 2^{m+1}}$.
\item[$\bullet$] $A = \bigcup \{\Omega_{j_1 \cdots j_m}: \, j_1, \dots ,j_m \in \{0,1\} \}$.
\end{enumerate}
As usual, we denote by \,$\chi_S$ \,the characteristic function of a set \,$S \subset \Omega$, that is,
$\chi_S(x) = \left\{\begin{array}{ll} 1 & \mbox{if }  x \in S, \\
                                       0 & \mbox{if }  x \not\in S.
                                       \end{array} \right.$
Let \,$n \in \N$ \,with \,$n \ge 2$. Then there is a unique \,$m = m(n) \in \N$ \,as well as unique \,$j_1, \dots ,j_m \in \{0,1\}$ \,such that
$$
n = 2^m + 2^{m-1}j_1 + 2^{m-2}j_2 + \cdots + 2j_{m-1} + j_m.
$$
Define \,$f_1 := \chi_A$ \,and \,$f_n := \chi_{\Omega_{j_1 \cdots j_m}}$.
Observe that \,$m(n) \to \infty$ \,as \,$n \to \infty$. For fixed \,$q \in (0,+\infty )$ \,we have
$$
\lim_{n \to \infty} \int_\Omega |f_n|^q \, d\mu = \lim_{n \to \infty} \mu (\Omega_{j_1 \cdots j_m}) = \lim_{n \to \infty} {\mu (A) \over 2^{m(n)}} = 0,
$$
which yields \,$(f_n) \in \bigcap_{q > 0} \mathcal{S}_{L_q}$. Now, fix \,$x_0 \in A$ \,and \,$m \in \N$.
Then, there is a unique \,$(j_1, \dots, j_m) \in \{0,1\}^m$ \,such that \,$x_0 \in \Omega_{j_1 \cdots j_m}$.
If \,$0 \le k < 2^m$ \,then
$$
f_{2^m + k}(x_0) = \left\{\begin{array}{ll} 1 & \mbox{if }  k = 2^{m-1}j_1 + 2^{m-2}j_2 + \cdots + 2j_{m-1} + j_m, \\
                                       0 & \mbox{otherwise.}
                                       \end{array} \right.
$$
Consequently, the sequence \,$(f_n(x_0))$ \,contains infinitely many $0$'s and infinitely many $1$'s, and so it does not converge.
This tells us that \,$(f_n) \not\in \mathcal{S}_p$. It follows that \,$\Big( \bigcap_{q > 0} \mathcal{S}_{L_q} \Big) \setminus \mathcal{S}_{p} \ne \varnothing$.
The $\mathfrak{c}$-lineability of \,$\Big( \bigcap_{q > 0} \mathcal{S}_{L_q} \Big) \setminus \mathcal{S}_{p}$ \,follows after applying
Theorem \ref{Tma de Vecina dense-lineable} on \,$V = L_0$, $\mathcal{C} = \bigcap_{q > 0} \mathcal{S}_{L_q}$ \,and \,$\mathcal{D} = \mathcal{S}_{p}$.

\vskip 3pt

Finally, our assumptions together with Lemma \ref{Lemma L0N metrizable and separable} yield that \,$L_0^{\N}$ \,is metrizable and separable.
Thus, the desired $\mathfrak{c}$-dense lineability of our family is directly deduced from Theorem \ref{A stronger than B} just by applying it on
\,$X = L_0^{\N}$, $A = \Big( \bigcap_{q > 0} \mathcal{S}_{L_q} \Big) \setminus \mathcal{S}_{p}$ \,and \,$B = c_{00}(L_0)$.

\vskip 3pt

The $\mathfrak{c}$-dense lineability of the sets of the last assertion of the theorem comes from the facts that
\,$\mathcal{S}_c \subset \mathcal{S}_{p}$,
$\mathcal{S}_{u} \subset \mathcal{S}_{au}  \subset \mathcal{S}_{p}$, and that any
\,$\mathcal{S}_{L_q}$ \,is contained in \,$\mathcal{S}_m$.
\end{proof}

Two new auxiliary results will be used in the next theorem. For each \,$A \in \mathcal{A}$, we can consider the induced measure space
\,$(A,\mathcal{A}_A,\mu |_{\mathcal{A}_A})$ \,as well as the corresponding space \,$L_0(A)$ \,of measurable functions \,$A \to \R$,
endowed with the local convergence in measure. Here \,$\mathcal{A}_A$ \,stands for the induced $\sigma$-algebra \,$\{M \in \mathcal{A}: \, M \subset A \}$.
Moreover, we define the vector subspace
$$
L_{0,A} := \{ f \in L_0 : \, f(x) = 0 \hbox{ for almost every } x \in \Omega \setminus A \}.
$$
From the fact that convergence in measure of a sequence implies pointwise a.e.~convergence of a subsequence, it follows that
\,$L_{0,A}$ \,is closed in \,$L_0$. In addition, it is easy to see that the mapping \,$f \in L_0(A) \mapsto f|_A \in L_{0,A}$,
where the last set is endowed with the inherited topology from $L_0$,
is an algebraic and topological isomorphism from \,$L_0(A)$ \,onto $L_{0,A}$. Consequently, we obtain the following.

\begin{lemma} \label{L_0 (A)^N}
Let \,$A \in \mathcal{A}$. Then we have:
\begin{enumerate}
\item[\rm (a)] $L_{0,A}^{\N}$ \,is a closed subset of \,$L_0^{\N}$.
\item[\rm (b)] The topological vector spaces \,$L_0(A)^{\N}$ \,and \,$L_{0,A}^{\N}$ \,are isomorphic, hence homeomorphic.
In particular, if \,$A$ \,is $\sigma$-finite, then \,$L_{0,A}^{\N}$ \,is metrizable.
\end{enumerate}
\end{lemma}

Recall that the {\it support} of a function \,$\varphi : \Omega \to \R$ \,is the set
$$\sigma (\varphi ) := \{x \in \Omega : \, \varphi (x) \ne 0 \}.$$

\begin{lemma} \label{Lemma closure in L0 de span disjoint}
Assume that \,$(\Omega, \mathcal{A}, \mu )$ is a measure space and that \,$(\varphi_n) \in L_0^{\N}$ \,is a sequence of measurable functions whose supports
\,$\sigma (\varphi_n)$ \,are mutually disjoint and have finite measure. Then
$$
\overline{\rm span} \{\varphi_n : \, n \in \N \} = \Big\{ \sum_{n=1}^{\infty} c_n \varphi_n : \, c_n \in \R \hbox{ \,for all } \,n \in \N \Big\},
$$
where the closure is in \,$L_0$ \,under the topology of local convergence in measure.
\end{lemma}

\begin{proof}
Notice that without loss of generality we can assume that \,$\mu (\sigma (\varphi_n)) > 0$ \,for all \,$n \in \N$.
Let \,$A := \bigcup_{n=1}^\infty \sigma (\varphi_n)$. Then \,$A \in \mathcal{A}$ \,and is $\sigma$-finite.
Hence \,$L_{0,A}$ \,is metrizable.
Denote \,$\mathcal{S} := \Big\{ \sum_{n=1}^{\infty} c_n \varphi_n : \, c_n \in \R \hbox{ \,for all } \,n \in \N \Big\}$.
Then, trivially, the sets
\,${\rm span} \{\varphi_n : \, n \in \N \}$ \,and \,$\mathcal{S}$ \,are contained in \,$L_{0,A}$.
Since the last set is closed in \,$L_0$, the closure \,$\mathcal{S}_0$ \,of \,${\rm span} \{\varphi_n : \, n \in \N \}$ \,in \,$L_{0,A}$
\,equals its closure in \,$L_0$. Then, it should be shown that \,$\mathcal{S}_0 = \mathcal{S}$.

\vskip 3pt

Let us prove the required equality by double inclusion. Let \,$F = \sum_{n=1}^{\infty} c_n \varphi_n \in \mathcal{S}$
\,and set \,$F_n := \sum_{k=1}^{n} c_k \varphi_k \in {\rm span} \{\varphi_j : \, j \in \N \}$ \,($n \in \N$).
We have that \,$F_n \to F$ \,locally in measure (which implies \,$\mathcal{S} \subset \mathcal{S}_0$). Indeed,
fix \,$K \in \mathcal{A}$ \,with \,$\mu (K) < \infty$, and define \,$K_n := K \cap (\sigma (\varphi_1) \cup \sigma (\varphi_2) \cup \cdots \cup \sigma (\varphi_n))$
\,and \,$\widetilde{K} := K \cap \bigcup_{n=1}^\infty \sigma (\varphi_n)$.
Observe that the sequence \,$(K_n)$ \,increases to \,$\widetilde{K}$. Therefore, $\mu (K_n) \to \mu (\widetilde{K})$, and so \,$\mu (\widetilde{K} \setminus K_n) \to 0$.
Let \,$\ve > 0$.
Since \,$F$ \,and the \,$F_n$'s \,vanish outside \,$\bigcup_{n=1}^\infty \sigma (\varphi_n)$, we get
\begin{equation*}
\begin{split}
\{x \in K : \, |F_n(x) - F(x)| > \ve \}& = \{x \in \widetilde{K} : \, |F_n(x) - F(x)| > \ve \} \\
                                                   &\subset \{x \in \widetilde{K} : \, F_n(x) \ne F(x) \} \subset \widetilde{K} \setminus K_n.
\end{split}
\end{equation*}
It follows that
$$
\mu (\{x \in K : \, |F_n(x) - F(x)| > \ve \}) \longrightarrow 0 \hbox{ \,as } \,n \to \infty,
$$
that is, $F_n \to F$ \,locally in measure.

\vskip 3pt

As for the reverse inclusion, suppose that \,$F \in \mathcal{S}_0$. Since \,$L_{0,A}$ \,is metrizable, there exists a sequence
\,$(F_n) \subset {\rm span} \{\varphi_n : \, n \in \N \}$ \,such that \,$F_n \to F$ \,in measure.
For each \,$n \in \N$ \,there exists \,$\nu (n) \in \N$ \,as well as scalars \,$c_{n,j}$ \,($j=1, \dots ,\nu (n)$)
\,such that \,$F_n = \sum_{j=1}^{\nu (n)} c_{n,j} \varphi_j$. If we define \,$c_{n,j} := 0$ \,for \,$j > \nu (n)$, then we can write
\,$F_n = \sum_{j=1}^{\infty} c_{n,j} \varphi_j$. Due to convergence in measure, we can select a sequence \,$n(1) < n(2) < \cdots < n(k) < \cdots$
\,of natural numbers as well as a $\mu$-null set \,$Z \in \mathcal A$ \,such that
\begin{equation} \label{Equation Ccia puntual de una subs}
F_{n(k)}(x) \longrightarrow F(x) \,\hbox{ as } \,k \to \infty, \hbox{\,for every } x \in \Omega \setminus Z.
\end{equation}
Now, fix \,$l \in \N$ \,and choose \,$x_l \in \sigma (\varphi_l) \setminus Z$. It follows from \eqref{Equation Ccia puntual de una subs} that
$$
c_{n(k),l} \cdot \varphi_{l}(x_l) = F_{n(k)}(x_l) \longrightarrow F(x_l) \hbox{ \,as} \,k \to \infty.
$$
Since \,$\varphi_{n(k)}(x_l) \ne 0$, we get \,$\lim_{k \to \infty} c_{n(k),l} = {F(x_l) \over \varphi_{l}(x_l)}$.
But the \,$c_{n(k),l}$'s do not depend on the chosen point $x_l \in \sigma (\varphi_l) \setminus Z$.
Consequently, there is \,$c_l \in \R$ \,such that \,$F(x) = c_l \cdot \varphi_l(x)$ \,for almost every \,$x \in \sigma (\varphi_k)$.
Since the supports of the \,$\varphi_n$'s are pairwise disjoint, this entails that \,$F = \sum_{n=1}^{\infty} c_n \varphi_n$.
In other words, $F \in \mathcal{S}$, as desired.
\end{proof}

\begin{theorem} \label{Sm-Spae spaceable}
Let \,$(\Omega ,\mathcal{A},\mu )$ \,be a nonatomic semifinite measure space.
Then \,$\mathcal{S}_m \setminus \mathcal{S}_{p}$ \,is spaceable in \,$L_0^{\N}$.
In particular, the sets \,$\mathcal{S}_m \setminus \mathcal{S}_{c}$, $\mathcal{S}_m \setminus \mathcal{S}_{au}$ \,and \,$\mathcal{S}_m \setminus \mathcal{S}_{u}$
\,are spaceable in \,$L_0^{\N}$.
\end{theorem}

\begin{proof}
The last assertion follows from the inclusions \,$\mathcal{S}_{u} \subset \mathcal{S}_{au} \subset \mathcal{S}_{p}$ \,and $\mathcal{S}_c \subset \mathcal{S}_{p}$.
Thus, our task is to show the spaceability of \,$\mathcal{S}_m \setminus \mathcal{S}_{p}$.

\vskip 3pt

Our construction is a refinement of the one given in the proof of Theorem \ref{Sm-Spae dense lineable}.
By Proposition \ref{Prop-caracterization semifinite, and atomless}, there are sets
\,$A,S \in \mathcal{A}$ \,such that \,$0 < \mu (A) < \infty$, $S \subset A$, and \,$\mu(S) = {\mu (A) \over 2}$.
Let \,$A_1 := S$. Since \,$\mu (A \setminus S) = {\mu (A) \over 2}$, again by \ref{Prop-caracterization semifinite, and atomless}
we can find a measurable set \,$A_2 \subset A \setminus S = A \setminus A_1$ (hence $A_1 \cap A_2 = \varnothing$) \,such that \,$\mu(A_2) = {\mu (A) \over 2^2}$.
Now, we select \,$A_3 \in \mathcal{A}$ \,with \,$A_3 \subset A \setminus (A_1 \cup A_2)$ \,and \,$\mu(A_3) = {\mu (A) \over 2^3}$.
Continuing this process, we obtain the existence of a countable family \,$A_1,A_2, \dots,A_l, \dots$ \,of mutually disjoint measurable sets such
that \,$\mu (A_l) = {\mu (A) \over 2^l}$ \,for all \,$l \in \N$. Now, we apply to each \,$A_l$ \,the decomposition method given in the proof of
the preceding theorem. Therefore, we obtain for every \,$l \in \N$ \,a family
$$
\big\{ \Omega_{lj_1 \cdots j_m} : \, m \in \N, \, j_1, \dots ,j_m \in \{0,1\} \big\} \subset \mathcal{A}
$$
satisfying the following properties for each \,$m \in \N$:
\begin{enumerate}
\item[$\bullet$] $\Omega_{lj_1 \cdots j_m} = \Omega_{lj_1 \cdots j_m 0} \cup \Omega_{lj_1 \cdots j_m 1}$.
\item[$\bullet$] $\Omega_{lj_1 \cdots j_m 0} \cap \Omega_{lj_1 \cdots j_m 1} = \varnothing$.
\item[$\bullet$] $\mu (\Omega_{lj_1 \cdots j_m 0}) =  \mu(\Omega_{lj_1 \cdots j_m1}) = {\mu (A) \over 2^{m+1+l}}$.
\item[$\bullet$] $A_l = \bigcup \{\Omega_{lj_1 \cdots j_m}: \, j_1, \dots ,j_m \in \{0,1\} \}$.
\end{enumerate}
We are going to define a countable set \,$\{ {\bf f}_l : \,l \in \N \}$ \,of sequences \,${\bf f}_l = (f_{l,n}) \in L_0^{\N}$.
To this end, let us set \,$f_{l,1} := \chi_{A_l}$. Now, fix \,$l \in \N$ \,and \,$n \ge 2$.
Then, there is a unique \,$m = m(n) \in \N$ \,as well as unique \,$j_1, \dots ,j_m \in \{0,1\}$ \,such that
\,$n = 2^m + 2^{m-1}j_1 + 2^{m-2}j_2 + \cdots + 2j_{m-1} + j_m$. Let us set \,$f_{l,n} := \chi_{\Omega_{lj_1 \cdots j_m}}$.
This way, the sequences \,${\bf f}_l$ \,have been defined.

\vskip 3pt

Let us define
$$
M := \overline{\rm span} \{ {\bf f}_l : \, l \in \N \},
$$
the closure in \,$L_0^{\N}$ \,of the linear span of the sequences \,${\bf f}_l$. Plainly, $M$ is a closed vector subspace of \,$L_0^{\N}$.
From the fact that the sets \,$A_l$ \,are mutually disjoint, it is easily derived the linear independence of the first coordinates \,$f_{1,1},f_{2,1},\dots ,f_{3,1}, \dots$,
which in turn implies the linear independence of the family \,$\{{\bf f}_l : \, l \in \N \}$. Hence \,$M$ \,is infinite dimensional.
Thus, our unique task is to prove that \,$M \setminus \{0\} \subset \mathcal{S}_{m} \setminus \mathcal{S}_{p}$.

\vskip 3pt

To this end, fix \,${\bf F} = (F_n) \in M \setminus \{0\}$. Observe first that if we denote \,$C := \bigcup_{l \in \N} A_l$ \,then \,$C$ \,is $\sigma$-finite
and all sequences \,${\bf f}_l$ \,are in \,$L_{0,C}^{\N}$, so ${\rm span} \{ {\bf f}_l : \, l \in \N \} \subset L_{0,C}$.
By Lemma \ref{L_0 (A)^N}, $L_{0,C}^{\N}$ \,is metrizable and \,$M \subset L_{0,C}^{\N}$.
Consequently, there exists a sequence \,$({\bf g}_l) \subset {\rm span} \{ {\bf f}_n : \, n \in \N \}$
\,such that \,${\bf g}_l \to {\bf F}$ \,in \,$L_0^{\N}$. Letting \,${\bf g}_l = (g_{l,n})$,
we have for each \,$n \in \N$ \,that \,$g_{l,n} \to F_n$ \,in measure. This implies that \,$F_n \in \overline{\rm span} \{f_{l,n}: \, l \in \N \}$.
Thanks to Lemma \ref{Lemma closure in L0 de span disjoint}, we obtain for every \,$n \in \N$ \,the existence of a real sequence
\,$(\alpha_{l,n})$ \,such that \,$F_n = \sum_{l=1}^{\infty} \alpha_{l,n} f_{l,n}$.

\vskip 3pt

On the one hand, we have for each \,$\varepsilon > 0$ \,that
\begin{equation*}
\begin{split}
\mu (|F_n - 0| > \varepsilon ) &\le \mu \big( \bigcup_{l=1}^\infty \{ f_{l,n} \ne 0 \} \big) = \sum_{l=1}^{\infty} \mu (\Omega_{lj_1 \cdots j_m}) \\
                               &= \sum_{l=1}^{\infty} {\mu (A) \over 2^{m(n)+l}} = \frac{\mu (A)}{2^{m(n)}} \longrightarrow 0 \hbox{ \,as } \,n \to \infty.
\end{split}
\end{equation*}
This shows that \,$F_n \to 0$ \,in measure, that is, ${\bf F} \in \mathcal{S}_{m}$. On the other hand, a closer look to the proof of
Lemma \ref{Lemma closure in L0 de span disjoint} together with the fact that the coordinates \,$g_{l,n}$ ($n=1,2, \dots$) \,inherit the same coefficients of each \,${\bf g}_l$
\,(as a linear combination of ${\bf f}_1,{\bf f}_2, \dots$) reveals that the coefficients \,$\alpha_{l,n}$ \,are the same for all \,$F_n$, that is, there is a real sequence
\,$(\alpha_l)$ \,with
$$
F_n = \sum_{l=1}^{\infty} \alpha_{l} f_{l,n} \ \ \hbox{for all } \, n \in \N .
$$
Since \,${\bf F} \ne 0$, there is \,$L \in \N$
\,such that \,$\alpha_L \ne 0$. Recall that \,$f_{L,n} := \chi_{\Omega_{Lj_1 \cdots j_m}}$, with \,$m$ \,and \,$j_1, \dots ,j_m$ \,depending on \,$n$.
If \,$x_0 \in A_L$, the sequence
$$
F_n(x_0) = \alpha_L \cdot \chi_{\Omega_{Lj_1 \cdots j_m}} (x_0) \quad (n=1,2,3, \dots )
$$
contains infinitely many zeros and infinitely many \,$\alpha_L$'s.
Since \,$\mu (A_L) > 0$, the sequence \,$(F_n)$ \,does not tend to \,$0$ \,pointwise a.e. In other words, ${\bf F} \not\in \mathcal{S}_{p}$, and so
${\bf F} \in \mathcal{S}_m  \setminus \mathcal{S}_{p}$, as required.
\end{proof}

It is possible to reinforce the lineability properties of \,$\mathcal{S}_m \setminus \mathcal{S}_p$ \,that have been furnished in the preceding two theorems, so as to exclude not
only a.e.~pointwise convergence but also any $q$-mean convergence.

\begin{theorem} \label{Sm-(Sp u USLq) dense lineable y spaceable}
Let \,$(\Omega ,\mathcal{A},\mu )$ \,be a nonatomic measure space.
\begin{enumerate}
\item[\rm (a)] If the measure space is $\sigma$-finite and satisfies {\rm (S)}, then
$\mathcal{S}_m \setminus \Big( \mathcal{S}_p \cup \bigcup_{q > 0} \mathcal{S}_{L_q} \Big)$
is $\mathfrak{c}$-dense-lineable in \,$L_0^{\N}$.
\item[\rm (b)] If the measure space is semifinite, then $\mathcal{S}_m \setminus \Big( \mathcal{S}_p \cup \bigcup_{q > 0} \mathcal{S}_{L_q} \Big)$
is spaceable in $L_0^{\N}$.
\end{enumerate}
\end{theorem}

\begin{proof}
To prove (a), we shall follow the procedure (and keep the notation) of the proof of Theorem \ref{Sm-Spae dense lineable}, so that we start with
a measurable set \,$A$ \,with \,$0 < \mu (A) < \infty$, and we define \,${\bf f} = (f_n) \in L_0^{\N}$ \,as
\,$f_1 = \chi_A$ \,and \,$f_n = \chi_{\Omega_{j_1 \cdots j_m}}$ ($n \ge 2$), where as \,$m$ \,as the \,$j_i$'s \,depend on \,$n$.
As in the mentioned theorem, it suffices, in order to prove (a), to exhibit a
sequence
$$
{\bf g} \in \mathcal{S}_m \setminus \Big( \mathcal{S}_p \cup \bigcup_{q > 0} \mathcal{S}_{L_q} \Big) .
$$
Indeed, after having found such a sequence, the application of Theorems \ref{A stronger than B} and \ref{Tma de Vecina dense-lineable}
would be similar: just take into account that \,$\mathcal{S}_p \cup \bigcup_{q > 0} \mathcal{S}_{L_q}$ \,satisfies conditions (a) to (e)
in Theorem \ref{Tma de Vecina dense-lineable}
\,and that \,$c_0(L_0) + \big( \mathcal{S}_p \cup \bigcup_{q > 0} \mathcal{S}_{L_q} \big) \subset \mathcal{S}_p \cup \bigcup_{q > 0} \mathcal{S}_{L_q}$.
Now, hands to work. If we define \,${\bf g} = (g_n)$ \,as
$$
g_n = 2^{m(n)^2} \cdot f_n,
$$
then we have:
\begin{itemize}
\item[$\bullet$] Given \,$\ve > 0$, it is evident that \,$\{ |g_n| > \ve \} \subset \{ f_n \ne 0 \} = \Omega_{j_1 \cdots j_m}$, which implies
$$
\mu (\{ |g_n - 0| > \ve \}) \le  \mu ( \Omega_{j_1 \cdots j_m} ) = {\mu (A) \over 2^{m(n)}} \longrightarrow 0 \hbox{ \ as } \, n \to \infty
$$
because \,$m(n) \to \infty$. Then ${\bf g} \in \mathcal{S}_m$.
\item[$\bullet$] Given \,$x_0 \in A$, we can see as in the proof of Theorem \ref{Sm-Spae dense lineable} that there are infinitely many
\,$n \in \N$ \,such that \,$f_n(x_0) = 1$, and so \,$g_n(x_0) = 2^{m(n)^2}$ \,for such $n$'s. Then \,$(g_n(x_0))$ \,cannot tend to zero.
Since \,$\mu (A) > 0$, we get \,${\bf g} \not \in \mathcal{S}_p$.
\item[$\bullet$] Given \,$q > 0$, we have that
$$
\lim_{n \to \infty} \int_{\Omega} |g_n|^q \, d\mu = \lim_{n \to \infty} 2^{m(n)^2 q} \cdot  \mu ( \Omega_{j_1 \cdots j_m} )
= \lim_{n \to \infty} 2^{m(n)^2 q} \cdot  {\mu (A) \over 2^{m(n)}} = \infty .
$$
\end{itemize}
Therefore, ${\bf g} \not\in \mathcal{S}_{L_q}$ \,and, thus, ${\bf g} \in \mathcal{S}_m \setminus \Big( \mathcal{S}_p \cup \bigcup_{q > 0} \mathcal{S}_{L_q} \Big)$, as required.

\vskip 3pt

In order to prove (b), we shall follow the approach (and keep the notation) of the proof of Theorem \ref{Sm-Spae spaceable}, so that
we start with a collection
\,$\{A_l : \, n \in \N \} \subset \mathcal{A}$ \,satisfying \,$\mu (A_l) = {\mu (A) \over 2^l}$, where \,$0 < \mu (A) < \infty$.
As in that proof, we consider the sequences \,${\bf f}_l = (f_{l,n}) \in L_0^{\N}$ ($l=1,2, \dots $)
\,given by
$$
f_{l,1} = \chi_{A_l} \ \ {\rm and} \ \ f_{l,n} = \chi_{\Omega_{lj_1 \cdots j_m}} \ (n \ge 2).
$$
Now, define \,${\bf g}_l = (g_{l,n})$ ($l = 1,2, \dots $) \,as \,$g_{l,n} = 2^{m(n)^2}f_{l,n}$ ($n=1,2, \dots $), and let
$$
M := \overline{\rm span} \, \{ {\bf g}_l : \, l \in \N \}.
$$
Note that, as in the mentioned proof, $M$ is a closed infinite dimensional vector subspace of \,$L_0^{\N}$ \,and, given \,${\bf G} = (G_n) \in M \setminus \{0\}$,
there exists a real sequence \,$(\alpha_l)$ \,with \,$\alpha_L \ne 0$ \,for some \,$L \in \N$ \,and
\,$G_n = \sum_{l=1}^{\infty} \alpha_{l} g_{l,n} = 2^{m(n)^2} \cdot \sum_{l=1}^{\infty} \alpha_{l} f_{l,n}$
\,for all \,$n \in \N$.

\vskip 3pt

On the one hand, we have for each \,$\varepsilon > 0$ \,that
\begin{equation*}
\begin{split}
\mu (|G_n - 0| > \varepsilon ) &\le \mu \big( \bigcup_{l=1}^\infty \{ g_{l,n} \ne 0 \} \big) = \sum_{l=1}^{\infty} \mu (\Omega_{lj_1 \cdots j_m}) \\
                               &= \sum_{l=1}^{\infty} {\mu (A) \over 2^{m(n)+l}} = \frac{\mu (A)}{2^{m(n)}} \longrightarrow 0 \hbox{ \,as } \,n \to \infty ,
\end{split}
\end{equation*}
which shows that \,${\bf G} \in \mathcal{S}_{m}$.
On the other hand, if \,$x_0 \in A_L$, the sequence \,$(G_n(x_0)) = \big( \alpha_L \cdot 2^{m(n)^2} \cdot \chi_{\Omega_{Lj_1 \cdots j_{m(n)}}} \big)$
\,contains infinitely many terms \,$2^{m(n)^2} \alpha_L$.
Since \,$\mu (A_L) > 0$, the sequence \,$(G_n)$ \,does not tend to \,$0$ \,pointwise almost everywhere.
In other words, ${\bf G} \not\in \mathcal{S}_{p}$.

\vskip 3pt

Finally, given \,$q > 0$, we have that
\begin{equation*}
\begin{split}
\int_{\Omega} |G_n|^q \, d\mu &\ge \int_{A_L} |G_n|^q \, d\mu =  \int_{\Omega_{Lj_1 \cdots j_{m(n)}}} 1^q \,d\mu \\
                              &= 2^{m(n)^2 q} \cdot \mu (\Omega_{Lj_1 \cdots j_{m(n)}}) = 2^{m(n)^2 q - L - m(n)} \longrightarrow \infty \hbox{ \ as } \, n \to \infty .
\end{split}
\end{equation*}
Thus, ${\bf G} \not \in \mathcal{S}_{L_q}$.
Consequently, $M \setminus \{0\} \subset \mathcal{S}_m \setminus \Big( \mathcal{S}_p \cup \bigcup_{q > 0} \mathcal{S}_{L_q} \Big)$, which
yields the spaceability of the last set.
\end{proof}

Observe that, in order to obtain lineability for \,$\mathcal{S}_{p} \setminus \mathcal{S}_m$, the measure must be infinite, because \,$\mathcal{S}_{au} \subset \mathcal{S}_m$ and, if the measure is finite, then \,$\mathcal{S}_{au} = \mathcal{S}_p$ \,due to Egoroff's theorem. We are going to consider the following two properties on a measure space \,$(\Omega , \mathcal{A}, \mu )$, both of them implying \,$\mu (\Omega ) = \infty$:

\vskip 3pt

\noindent \phantom{aaa} (P) {\it There is a countable family \,$\{A_n : \, n \in \N \} \subset \mathcal A$ \,of mutually \phantom{aaaaaa} \break
\phantom{aaaaaaa\,} disjoint measurable sets such that \,$\inf_{n \in \N} \mu (A_n) > 0$.}

\vskip 3pt

\noindent \phantom{aaa} (Q) \ $\sup \{\mu (S) : \,S \in \mathcal{A}$ such that $\mu (S) < \infty \} = \infty$.

\vskip 3pt

On the one hand, it follows from \cite[Theorem 14.23]{nielsen} that (Q) holds if and only if there is a countable family \,$\{A_n : \, n \in \N \} \subset \mathcal A$
\,of mutually disjoint measurable sets such that \,$1 \le \mu (A_n) < \infty$ \,for all \,$n \in \N$. Therefore, (Q) implies (P) (the converse is not true: consider any infinite set \,$\Omega$ \,and define on \,$\mathcal P(\Omega )$ \,the measure \,$\mu (S) = \infty$ \,for all \,$S \subset \Omega$ \,with \,$S \ne \varnothing$). On the other hand,
from Lemma \ref{Prop-caracterization semifinite, and atomless} it is easy to extract that any nonatomic semifinite measure space with \,$\mu (\Omega ) = \infty$ \,satisfies (Q), hence (P). But being nonatomic is not at all necessary for (P) or (Q): consider, for instance, the counting measure \,$\mu$ \,on the measurable space \,$(\N, \mathcal{P} (\N ))$ \,and select the sets \,$A_n = \{n\}$ \,($n \ge 1$). Observe that, in addition, that \,$(\N, \mathcal{P} (\N ), \mu )$ \,is $\sigma$-finite and satisfies condition (S).

\begin{theorem} \label{Spae-Sm dense lineable}
Let \,$(\Omega ,\mathcal{A},\mu )$ \,be a $\sigma$-finite measure space satisfying conditions {\em (S)} and {\em (P).}
Then \,$\mathcal{S}_{p} \setminus \mathcal{S}_m$ \,is $\mathfrak{c}$-dense-lineable in \,$L_0^{\N}$.
In particular, $\mathcal{S}_{p} \setminus \bigcup_{q > 0} \mathcal{S}_{L_q}$, $\mathcal{S}_{p} \setminus \mathcal{S}_{c}$,
$\mathcal{S}_{p} \setminus \mathcal{S}_{au}$ \,and \,$\mathcal{S}_{p} \setminus \mathcal{S}_u$ \,are $\mathfrak{c}$-dense-lineable in \,$L_0^{\N}$.
\end{theorem}

\begin{proof}
Define the sequence \,${\bf f} = (f_n) \in L_0^{\N}$ \,as \,$f_n := \chi_{A_n}$, where \,$(A_n)$ \,is the sequence of measurable sets furnished by condition (P).
Then there exists \,$\alpha > 0$ \,such that \,$\mu (A_n) \ge \alpha$ \,for all \,$n \in \N$.

\vskip 3pt

Let \,$x_0 \in \Omega$. Since the \,$A_n$'s \,are pairwise disjoint, $x_0$ belongs to at most one \,$A_n$. Therefore, there is \,$n_0 \in \N$ \,such that
\,$f_n(x_0) = \chi_{A_n}(x_0) = 0$ \,for all \,$n > n_0$. Then \,$f_n(x_0) \longrightarrow 0$ \,as \,$n \to \infty$, and this holds for all \,$x_0 \in \N$.
Thus, ${\bf f} \in \mathcal{S}_{p}$. Now, let \,$\varepsilon = 1/2$, and fix \,$n \in \N$. Then
$$
\mu (|f_n - 0| > 1/2) = \mu (f_n = 1) = \mu (A_n) \ge \alpha ,
$$
and so \,$\mu (|f_n - 0| > 1/2) \not\longrightarrow 0$ \,as \,$n \to \infty$. This tells us that \,${\bf f} \not\in \mathcal{S}_m$ \,and, consequently,
$\mathcal{S}_{p} \setminus \mathcal{S}_m \ne \varnothing$. The $\mathfrak{c}$-lineability of \,$\mathcal{S}_{p} \setminus \mathcal{S}_{m}$ \,follows after applying
Theorem \ref{Tma de Vecina dense-lineable} on \,$V = L_0$, $\mathcal{C} = \mathcal{S}_{p}$ \,and \,$\mathcal{D} = \mathcal{S}_{m}$.

\vskip 3pt

Finally, our assumptions together with Lemma \ref{Lemma L0N metrizable and separable} yield that \,$L_0^{\N}$ \,is metrizable and separable.
Thus, the desired $\mathfrak{c}$-dense lineability of our family is directly deduced from Theorem \ref{A stronger than B} just by applying it on
\,$X = L_0^{\N}$, $A = \mathcal{S}_{p} \setminus \mathcal{S}_{m}$ \,and \,$B = c_{00}(L_0)$.
To conclude, the last assertion of the theorem follows from the inclusions \,$\mathcal{S}_u \subset \mathcal{S}_{au} \subset \mathcal{S}_m$,
$\mathcal{S}_{c} \subset \mathcal{S}_{m}$
\,and \,$\bigcup_{q > 0} \mathcal{S}_{L_q} \subset \mathcal{S}_m$.
\end{proof}

\begin{theorem} \label{Spae-Sm spaceable}
Let \,$(\Omega ,\mathcal{A},\mu )$ \,be a measure space satisfying condition {\em (Q).}
Then \,$\mathcal{S}_{p} \setminus \mathcal{S}_m$ \,is spaceable in \,$L_0^{\N}$.
In particular, $\mathcal{S}_{p} \setminus \bigcup_{q > 0} \mathcal{S}_{L_q}$,
$\mathcal{S}_{p} \setminus \mathcal{S}_{c}$,
$\mathcal{S}_{p} \setminus \mathcal{S}_{au}$ \,and \,$\mathcal{S}_{p} \setminus \mathcal{S}_u$ \,are spaceable in \,$L_0^{\N}$.
\end{theorem}

\begin{proof}
The last assertion follows from the inclusions \,$\mathcal{S}_u \subset \mathcal{S}_{au} \subset \mathcal{S}_m$,
$\mathcal{S}_{c} \subset \mathcal{S}_{m}$ \,and \,$\bigcup_{q > 0} \mathcal{S}_{L_q} \subset \mathcal{S}_m$.
Therefore, it suffices to prove the spaceability of \,$\mathcal{S}_{p} \setminus \mathcal{S}_m$.
The assumption of (Q) together with the use of any bijection \,$\N \times \N \to \N$ \,yields the existence
of a double sequence \,$(A_{k,n})$ \,of mutually disjoint measurable sets satisfying \,$1 \le \mu (A_{k,n} ) < \infty$ \,for all \,$k,n \in \N$.
For each \,$k \in \N$ \,we define the sequence \,${\bf f}_k = (f_{k,n})$ \,by
$$
f_{k,n} := \chi_{A_{n,k}}.
$$
Let us set
$$
M := \overline{\rm span} \{ {\bf f}_k : \, k \in \N \},
$$
A similar argument to the one provided in the proof of Theorem \ref{Sm-Spae spaceable}, we get that
\,$M$ \,is an infinite dimensional closed vector subspace of \,$L_0^{\N}$, and
so we have only to show that \,$M \setminus \{0\} \subset \mathcal{S}_{p} \setminus \mathcal{S}_{m}$.

\vskip 3pt

With this aim, fix \,${\bf F} = (F_n) \in M \setminus \{0\}$ \,and denote \,$A := \bigcup_{k,n \in \N} A_{k,n}$.
Then \,$A$ \,is $\sigma$-finite and all sequences \,${\bf f}_k$ \,are in \,$L_{0,A}^{\N}$, so
\,${\rm span} \{ {\bf f}_k : \, k \in \N \} \subset L_{0,A}^{\N}$.
By Lemma \ref{L_0 (A)^N}, $L_{0,A}^{\N}$ \,is metrizable and \,$M \subset L_{0,A}^{\N}$.
Thus, there is a sequence \,$({\bf g}_k) \subset {\rm span} \{ {\bf f}_l : \, l \in \N \}$
\,such that \,${\bf g}_k \to {\bf F}$ \,in \,$L_0^{\N}$.
Letting \,${\bf g}_k = (g_{k,n})$,
we have for each \,$n \in \N$ \,that \,$g_{k,n} \to F_n$ \,in measure. This implies that \,$F_n \in \overline{\rm span} \{f_{k,n}: \, k \in \N \}$.
An application of Lemma \ref{Lemma closure in L0 de span disjoint}
and the fact that the coordinates \,$g_{k,n}$ ($n=1,2, \dots$) \,inherit the same coefficients of each \,${\bf g}_k$
\,(as a linear combination of ${\bf f}_1,{\bf f}_2, \dots$) \,yield
the existence of a real sequence
\,$(\alpha_{k})$ \,such that \,$F_n = \sum_{k=1}^{\infty} \alpha_{k} f_{k,n} = \sum_{k=1}^{\infty} \alpha_{k} \chi_{A_{k,n}}$.
It remains to prove that \,${\bf F} \in \mathcal{S}_{p} \setminus \mathcal{S}_m$.

\vskip 3pt

On the one hand, given \,$x_0 \in \Omega$, there is at most one \,$(k,n) \in \N \times \N$ \,such that \,$x_0 \in A_{k,n}$.
Then there is \,$n_0 \in \N$ \,such that \,$x_0 \not\in A_{k,n}$ \,for all \,$(k,n) \in \N \times \{n_0+1,n_0+2, \dots \}$, and so
for \,$n > n_0$ \,we have
$$
F_n(x_0) = \sum_{k=1}^{\infty} \alpha_{k} \chi_{A_{k,n}}(x_0) = \sum_{k=1}^{\infty} \alpha_{k} \cdot 0 = 0.
$$
Consequently, $F_n(x_0) \to 0$, which entails \,${\bf F} \in \mathcal{S}_{p}$.

\vskip 3pt

On the other hand, since \,${\bf F} \ne 0$ \,we can infer the existence of an \,$l \in \N$ \,such that \,$\alpha_l \ne 0$.
Let \,$\varepsilon := {|\alpha_l| \over 2} > 0$.
For any \,$n \in \N$ \,we have that \,$F_n(x) = \alpha_l$ \,for all \,$x \in A_{l,n}$, and so
$$
\mu (\{ x \in \Omega : \, |F_n(x) - 0| > \varepsilon \}) \ge \mu (A_{l,n}) \ge 1.
$$
Hence \,$(F_n)$ \,cannot converge to \,$0$ \,in measure, that is, ${\bf F} \not\in \mathcal{S}_{m}$, as required.
\end{proof}

In order to compare \,$\mathcal{S}_{au}$ \,with \,$\mathcal{S}_u$, we need the following natural condition, that in some sense is dual of (P):

\vskip 3pt

\noindent \phantom{aaa} (R) {\it There is a countable family \,$\{A_n : \, n \in \N \} \subset \mathcal A$ \,of mutually disjoint \break
\phantom{aaaaaaa\,} measurable sets such that \,$\mu (A_n) > 0$ \,for all \,$n \in \N$ \,and \,$\lim_{n \to \infty} \mu (A_n) = 0$.}

\vskip 3pt

\noindent By using Theorem 14.22 in \cite{nielsen} it is easy to see that (R) holds if and only if \,$\inf \{\mu (A): \, A \in \mathcal{A}$ and $\mu (A) > 0 \} = 0$. The counting measure \,$\mu$ \,on the measurable space \,$(\N, \mathcal{P} (\N ))$ \,does not satisfy (R). Nevertheless, the measure \,$\nu$ \,on the same measurable space given by \,$\nu (A) = \sum_{n \in A} 2^{-n}$ ($A \subset \N$) \,does fulfill (R). By Proposition \ref{Prop-caracterization semifinite, and atomless}, this condition is also satisfied by any nonatomic semifinite measure space.

\begin{theorem} \label{Sau-Su dense lineable}
Let \,$(\Omega ,\mathcal{A},\mu )$ \,be a $\sigma$-finite measure space satisfying conditions {\em (S)} and {\em (R).}
Then \,$(\mathcal{S}_c \cap \mathcal{S}_{au}) \setminus \mathcal{S}_u$ \,is $\mathfrak{c}$-dense-lineable in \,$L_0^{\N}$.
In particular, $\mathcal{S}_c \setminus \mathcal{S}_u$ \,and
\,$\mathcal{S}_{au} \setminus \mathcal{S}_u$ \,are $\mathfrak{c}$-dense-lineable in \,$L_0^{\N}$.
\end{theorem}

\begin{proof}
Assume for a moment that we have proved that \,$(\mathcal{S}_c \cap \mathcal{S}_{au}) \setminus \mathcal{S}_u$ \,is not empty.
Then the $\mathfrak{c}$-lineability of \,$(\mathcal{S}_c \cap \mathcal{S}_{au}) \setminus \mathcal{S}_u$ \,would follow after applying
Theorem \ref{Tma de Vecina dense-lineable} on \,$V = L_0$, $\mathcal{C} = \mathcal{S}_c \cap \mathcal{S}_{au}$ \,and \,$\mathcal{D} = \mathcal{S}_{u}$.
Moreover, thanks to Lemma \ref{Lemma L0N metrizable and separable} the space \,$L_0^{\N}$ \,is metrizable and separable.
Hence, the $\mathfrak{c}$-dense lineability of our family is deduced from Theorem \ref{A stronger than B} just by applying it on
the space \,$X = L_0^{\N}$ \,and the sets \,$A = (\mathcal{S}_c \cap \mathcal{S}_{au}) \setminus \mathcal{S}_u , \, B = c_{00}(L_0)$.

\vskip 3pt

Thus, it remains to prove that \,$(\mathcal{S}_c \cap \mathcal{S}_{au}) \setminus \mathcal{S}_u \ne \varnothing$.
To this end, the assumption yields the existence of countably many sets \,$A_n \in \mathcal{A}$ \,with \,$\mu (A_n) > 0$ \,for all \,$n \in \N$
\,and \,$\mu (A_n) \longrightarrow 0$ \,as \,$n \to \infty$. By extracting a subsequence if necessary, we may assume without loss of generality
that \,$\mu (A_n) < 2^{-n}$ \,for all \,$n \in \N$. Now we construct the following sequence of sets:
$$
B_n := \bigcup_{k=n}^\infty A_k \ \ (n=1,2,3, \dots ).
$$
Then \,$B_{n+1} \subset B_n$ \,and \,$0 < \mu (B_n) \le \sum_{k=n}^{\infty} \mu (A_k) < 2^{1-n} \longrightarrow 0$ \,as \,$n \to \infty$.
Define the sequence \,${\bf f} = (f_n) \in L_0^{\N}$ \,by
$$
f_n := \chi_{B_n} \ \ (n=1,2,3, \dots ).
$$
Let \,$\varepsilon > 0$. On the one hand, we can select \,$m \in \N$
\,with \,$2^{1-m} < \varepsilon$. Let \,$Z_\varepsilon := B_m \in \mathcal{A}$, that satisfies \,$\mu (Z_\varepsilon) < \ve$.
For every \,$x \in \Omega \setminus Z_\varepsilon$ \,we have that \,$x \not\in B_n$ \,for all \,$n \ge m$, and so \,$\chi_{B_n}(x) = 0$ \,for all \,$n \ge m$.
We infer that \,$\lim_{n \to \infty} \sup_{x \in \Omega \setminus Z_\varepsilon} |f_n(x) - 0| = 0$, that is, ${\bf f} \in \mathcal{S}_{au}$.
On the other hand, we have that \,$\{|f_n| > \ve \} \subset \{f_n \ne 0 \} = B_n$, and so
$$
\sum_{n=1}^{\infty} \mu (\{|f_n - 0| > \ve \}) \le \sum_{n=1}^{\infty} \mu (B_n) \le \sum_{n=1}^{\infty} 2^{1-n} = 2 < \infty .
$$
It follows that ${\bf f} \in \mathcal{S}_{c}$.

\vskip 3pt

Finally, let us suppose, via contradiction, that ${\bf f} \in \mathcal{S}_{u}$. Then there would exist a set \,$Z \in \mathcal{A}$
\,with \,$\mu (Z) = 0$ \,such that \,$f_n \to 0$ \,uniformly on \,$\Omega \setminus Z$. Observe that \,$(\Omega \setminus Z) \cap B_n \ne \varnothing$
\,for all \,$n \in \N$ (otherwise, $B_m \subset Z$ \,for some \,$m$, which is absurd because \,$\mu (Z) = 0$ \,and \,$\mu (B_m) > 0$).
Consequently, there is \,$x_n \in \Omega \setminus Z$ \,with \,$x_n \in B_n$. This entails \,$|f_n(x_n)| = 1$ ($n=1,2, \dots$),
which implies that \,$\sup_{x \in \Omega \setminus Z} |f_n(x)|$ \,cannot tend to $0$, which is absurd. Thus, ${\bf f} \not\in \mathcal{S}_{u}$
\,and this concludes the proof.
\end{proof}

\begin{theorem} \label{Sau-Su spaceable}
Let \,$(\Omega ,\mathcal{A},\mu )$ \,be a measure space satisfying condition {\em (R).} Then \break
$(\mathcal{S}_c \cap \mathcal{S}_{au}) \setminus \mathcal{S}_u$ \,is spaceable in \,$L_0^{\N}$.
In particular, $\mathcal{S}_c \setminus \mathcal{S}_u$ \,and
\,$\mathcal{S}_{au} \setminus \mathcal{S}_u$ \,are spaceable in \,$L_0^{\N}$.
\end{theorem}

\begin{proof}
By using the hypothesis (R) together with a bijection \,$\N \times \N \to \N$, we can obtain a family \,$A_{k,n}$ \,($k,n=1,2, \dots$) \,of mutually disjoint measurable sets satisfying \,$\mu (A_{k,n}) > 0$ \,for all \,$(k,n) \in \N \times \N$ \,and \,$\lim_{n \to \infty} \mu (A_{k,n}) = 0$ \,for all \,$k \in \N$.
In turn, by extracting appropriate subsequences, we may assume without loss of generality that \,$\mu (A_{k,n}) < 2^{-k-n-1}$ \,for all $(k,n) \in \N \times \N$.
Let us consider the following collection of measurable sets:
$$
B_{k,n} := \bigcup_{j=n}^\infty A_{k,j} \,\,\, (k,n=1,2, \dots ).
$$
Then \,$\mu (B_{k,n}) < 2^{-k-n}$. For each \,$k \in \N$, let \,${\bf f}_k = (f_{k,n}) \in L_0^{\N}$ \,be defined by
$$
f_{k,n} := \chi_{B_{k,n}}.
$$
The functions in the respective $1$st coordinates \,$f_{1,1},f_{2,1},f_{3,1}, \dots$ \,are linearly independent due to the disjointness of their supports
\,$A_{1,1},A_{2,1},A_{3,1}, \dots$. Hence the sequences \,${\bf f}_k$ \,are linearly independent, which entails that the set
$$
M := \overline{\rm span} \{ {\bf f}_k : \, k \in \N \}
$$
is a closed infinite dimensional vector subspace of \,$L_0^{\N}$. If we define
$$
A := \bigcup_{k,n=1}^{\infty} A_{k,n} = \bigcup_{k,n=1}^{\infty} B_{k,n},
$$
then \,${\bf f}_k \in L_{0,A}^{\N}$
\,and
$$
\mu (A) \le \sum_{k,n=1}^\infty \mu (B_{k,n}) < \sum_{k,n=1}^\infty 2^{-k-n} = 1 < \infty .
$$
By Lemma \ref{L_0 (A)^N}, the subspace \,$L_{0,A}^{\N}$
\,is a metrizable closed subspace of \,$L_{0,A}^{\N}$. Therefore, $M$ is the closure in \,$L_{0,A}^{\N}$ \,of \,${\rm span} \{ {\bf f}_k : \, k \in \N \}$.
To achieve the desired spaceability, it is enough to prove that \,$M \setminus \{0\} \subset \mathcal{S}_{au} \setminus \mathcal{S}_u$.
Let us fix \,${\bf F} = (F_n) \in M \setminus \{0\}$. Our goal is to show that \,${\bf F} \in \mathcal{S}_{au} \cap \mathcal{S}_{c}$ \,but \,${\bf F} \not\in \mathcal{S}_u$.
By Lemma \ref{Lemma closure in L0 de span disjoint} and arguing as in the proofs of Theorems \ref{Sm-Spae spaceable} and \ref{Spae-Sm spaceable},
we get the existence of reals \,$\alpha_1, \alpha_2, \dots , \alpha_k , \dots$ \,such that \,$F_n = \sum_{k=1}^{\infty} \alpha_k \chi_{B_{k,n}}$ \,for all \,$n \in \N$.
Since \,${\bf F} \ne 0$, there exists \,$l \in \N$ \,such that \,$\alpha_l \ne 0$.

\vskip 3pt

Given \,$\varepsilon > 0$, we can select \,$m \in \N$ \,such that \,$2^{-m+1} < \varepsilon$. Define
$$
Z_\varepsilon := \bigcup_{k=1}^\infty \bigcup_{n=m}^\infty B_{k,n} \in \mathcal{A}.
$$
Then
$$
\mu (Z_\varepsilon ) \le \sum_{k=1}^{\infty} \sum_{n=m}^{\infty} \mu (B_{k,n}) < \sum_{k=1}^{\infty} {2 \over 2^{k+m}} = {2 \over 2^m} < \varepsilon .
$$
Furthermore, for every \,$x \in \Omega \setminus Z_\varepsilon$ \,we have that \,$x \not \in B_{k,n}$ \,for all \,$(k,n) \in \N \times \N$
\,with \,$n \ge m$. Thus,
$$
\sup_{x \in Z_\varepsilon} |F_n(x) - 0| = \sup_{x \in Z_\varepsilon} \sum_{k=1}^{\infty} |\alpha_k| \, | \chi_{B_{k,n}}(x)| = 0
$$
for all \,$n \ge m$, and so \,$F_n \to 0$ \,uniformly on \,$\Omega \setminus Z_\varepsilon$, which proves that  \,${\bf F} \in \mathcal{S}_{au}$.
Now, observe that \,$\{|F_n| > \ve \} \subset \{F_n \ne 0 \} \subset \bigcup_{k \ge 1} B_{k,n}$, and so
$$
\sum_{n=1}^{\infty} \mu (\{|F_n - 0| > \ve \}) \le \sum_{n=1}^{\infty} \sum_{k=1}^{\infty} \mu (B_{k,n}) < \sum_{n=1}^{\infty} \sum_{k=1}^{\infty} 2^{-k-n} = 1 < \infty .
$$
From this we infer that ${\bf F} \in \mathcal{S}_{c}$.

\vskip 3pt

Finally, assume, by way of contradiction, that ${\bf F} \in \mathcal{S}_{u}$. This would imply the existence of
a measurable set \,$Z \in \mathcal{A}$ \,such that \,$\mu (Z) = 0$ \,and \,$F_n \to 0$ \,uniformly on \,$\Omega \setminus Z$.
Notice that \,$(\Omega \setminus Z) \cap B_{k,n} \ne \varnothing$
\,for all \,$k,n \in \N$ \,because, otherwise, $B_{k,n} \subset Z$ \,for some pair \,$(k,n)$, which is impossible since \,$\mu (Z) = 0$ \,and \,$\mu (B_{k,n}) > 0$.
In particular, we can choose for each \,$n \in \N$ \,a point \,$x_n \in \Omega \setminus Z$ \,with \,$x_n \in B_{l,n}$
(recall that $l \in \N$ \,has been chosen so as to satisfy $\alpha_l \ne 0$).
Since the sets \,$B_{k,n}$ \,are pairwise disjoint, it follows that \,$|F_n(x_n)| = |\alpha_l|$ ($n=1,2, \dots$),
which implies that \,$\sup_{x \in \Omega \setminus Z} |F_n(x)|$ \,cannot tend to $0$, which is absurd. Consequently, ${\bf F} \not\in \mathcal{S}_{u}$, as desired.
\end{proof}

To conclude this section, observe that when comparing uniform convergence to convergence in $q$-mean, it holds that if \,$\mu$ \,is finite then \,$\mathcal{S}_u \subset \mathcal{S}_{L_q}$ \,for all \,$q > 0$.
The point is that the finiteness of \,$\mu$ \,is, in an extreme way, the only condition for these inclusions.
Indeed, take the sequence of constant functions \,${\bf f} = (f_n) := (1/n)$. Then
\,${\bf f} \in  \mathcal{S}_u \setminus \bigcup_{q > 0} \mathcal{S}_{L_q} =: \mathcal{S}$, and an application of Theorem \ref{Tma de Vecina dense-lineable}
with \,$V = L_0^{\N}, \, \mathcal{C} = \mathcal{S}_u$ \,and \,$\mathcal{D} = \bigcup_{q > 0} \mathcal{S}_{L_q}$ \,yields the $\mathfrak{c}$-lineability of $\mathcal{S}$.
Finally, an application of Lemma \ref{Lemma L0N metrizable and separable} and of Theorem \ref{A stronger than B} with \,$X = L_0^{\N}$, $A = \mathcal{S}$ \,and
\,$B = c_{00}(L_0)$, together with the implications among the diverse kinds of convergence, gives the following result.

\begin{theorem} \label{Su-SL dense lineable}
Let \,$(\Omega ,\mathcal{A},\mu )$ \,be a $\sigma$-finite measure space satisfying condition {\em (S)} and $\mu (\Omega ) = \infty$.
Then \,$\mathcal{S}_u \setminus \bigcup_{q > 0} \mathcal{S}_{L_q}$ \,is $\mathfrak{c}$-dense-lineable in \,$L_0^{\N}$.
In particular, the sets \,$\mathcal{S}_p \setminus \bigcup_{q > 0} \mathcal{S}_{L_q}$, $\mathcal{S}_{m} \setminus \bigcup_{q > 0} \mathcal{S}_{L_q}$
and $\mathcal{S}_{au} \setminus \bigcup_{q > 0} \mathcal{S}_{L_q}$ \,are $\mathfrak{c}$-dense-lineable in \,$L_0^{\N}$.
\end{theorem}


\section{Algebrability}

\quad In this section it will be shown that, under not too strong constraints, the families of sequences we are considering contain, except for zero, large algebras.
In a natural way, the product operation that we shall consider in the algebra \,$L_0^{\N}$ \,is the one defined coordinatewise, that is, if \,${\bf f} = (f_n)$
\,and \,${\bf g} = (g_n)$ \,are in \,$L_0^{\N}$, then \,${\bf f} \cdot {\bf g} := (f_n g_n)$.

\vskip 3pt

As in the results of section 3, the final assertions of Theorems \ref{Sm-Spae algebrable}, \ref{Spae-Sm algebrable}
and \ref{Su-SL algebrable}
below follow from their respective main assertions together with the inclusions \,$\mathcal{S}_{L_q} \subset \mathcal{S}_m$ (for any $q > 0$),
$\mathcal{S}_c  \subset \mathcal{S}_m$, $\mathcal{S}_c \subset \mathcal{S}_p$,
$\mathcal{S}_{au} \subset \mathcal{S}_m$ \,and \,$\mathcal{S}_u  \subset \mathcal{S}_{au} \subset \mathcal{S}_p$.

\begin{theorem} \label{Sm-Spae algebrable}
Let \,$(\Omega ,\mathcal{A},\mu )$ \,be a nonatomic semifinite measure space. 
Then the set $\displaystyle{ \Big( \bigcap_{q > 0} \mathcal{S}_{L_q} \Big)} \setminus \mathcal{S}_{p}$ \,is strongly $\mathfrak{c}$-algebrable.
In particular, the sets \,$\mathcal{S}_m \setminus \mathcal{S}_{p}$,
$\mathcal{S}_m \setminus \mathcal{S}_{c}$,
$\mathcal{S}_m \setminus \mathcal{S}_{au}$ \,and \,$\mathcal{S}_m \setminus \mathcal{S}_{u}$
\,are strongly $\mathfrak{c}$-algebrable.
\end{theorem}

\begin{proof}
We have only to prove the strong $\mathfrak{c}$-algebrability of \,$\displaystyle{ \Big( \bigcap_{q > 0} \mathcal{S}_{L_q} \Big)} \setminus \mathcal{S}_{p}$.
Starting from the semifiniteness of \,$\mu$ \,and the absence of atoms, as in the proof of Theorem \ref{Sm-Spae dense lineable} (we keep the notation therein)
we can build a collection of measurable sets \,$\Omega_{j_1 \cdots j_m}$ ($j_1, \dots ,j_m \in \{0,1\}, \, m \in \N$) satisfying the properties given in
the mentioned proof. We had defined \,${\bf f} = (f_n) \in L_0^{\N}$ \,by \,$f_1 = \chi_A$ \,(recall that $0 < \mu (A) < \infty$) \,and \,$f_n = \chi_{\Omega_{j_1 \cdots j_m}}$ \,if \,$n \ge 2$,
where \,$m = m(n) \in \N$ \,and \,$j_1 = j_1(n), \dots , j_m = j_m (n)$ \,were uniquely determined by the binary decomposition \,$n = 2^m + 2^{m-1}j_1 + 2^{m-2}j_2 + \cdots + 2j_{m-1} + j_m$. Recall that \,$m(n) \to \infty$ \,as \,$n \to \infty$.

\vskip 3pt

The next construction is already classical in algebrability (see, for instance, algebrability proofs in \cite{polacos2023,calderongerlachpradomodes,vecina}).
Since the dimension of \,$\R$ \,as a vector space over the field \,$\Q$ \,of rationals is $\mathfrak{c}$, we can select a $\Q$-linearly independent set
\,$H \subset (0,+\infty )$ \,with ${\rm card} \, (H) = \mathfrak{c}$. Then the exponentials \,$t \mapsto e^{ct}$ ($c \in H$) \,will generate a free algebra.
We shall use a similar technique in the proofs of Theorems \ref{Spae-Sm algebrable} to \ref{Su-SL algebrable} below.

\vskip 3pt

In our case, we define the family
$$
\{ {\bf f}_c = (f_{c,n}) \in L_0^{\N} : \, c \in H \}
$$
by
$$
f_{c,n} = e^{c \cdot \sqrt{m(n)}} \cdot f_n.
$$
Denote \,$\N_0 := \N \cup \{0\}$. Assume that \,$N \in \N$ \,and that \,$P$ \,is a nonzero polynomial in \,$N$ \,variables without constant term.
Then, there is a nonempty finite set \,$J \subset \N_0^N \setminus \{(0,0, \dots ,0)\}$ \,and nonzero reals \,$a_{\bf j}$ (${\bf j} = (j_1, \dots ,j_N) \in J$)
\,such that
\begin{equation}\label{Equation Polynomial}
P(x_1, \dots ,x_N) = \sum_{{\bf j} \in J} a_{\bf j} \cdot x_1^{j_1} \cdots x_N^{j_N} \hbox{ \ for all } \,x_1, \dots ,x_N \in \R .
\end{equation}
Now, fix \,$N$ \,different elements of \,$H$, say \,$c_1,c_2, \dots ,c_N$. Since they are \,$\Q$-linearly independent, the numbers
\,$\sigma ({\bf j}) := c_1j_1 + \cdots + c_Nj_N$ (${\bf j} \in J$) \,are pairwise different. Then, there is one of them, say \,$\sigma ({\bf k})$,
such that \,$\sigma ({\bf k}) > \sigma ({\bf j})$ \,for all \,${\bf j} \in J \setminus \{{\bf k} \}$.
We have to prove that the sequence
$$
{\bf Q} = (Q_n) := P({\bf f}_{c_1}, \dots ,{\bf f}_{c_N})
$$
belongs to \,$\displaystyle{ \Big( \bigcap_{q > 0} \mathcal{S}_{L_q} \Big)} \setminus \mathcal{S}_{p}$ \,(that ${\bf Q} \ne 0$ is automatically implied by this).

\vskip 3pt

With this aim, note that for each \,$n \in \N$ \,we have
$$
Q_n = \sum_{{\bf j} \in J} a_{\bf j} \cdot f_{c_1,n}^{j_1} \cdots f_{c_N,n}^{j_N} =
\left( \sum_{{\bf j} \in J} a_{\bf j} \cdot e^{\sigma ({\bf j}) \sqrt{m(n)}}  \right) \cdot f_n.
$$
We have used the fact that each \,$f_n$ \,is the characteristic function of a set, and so \,$f_n^s = f_n$ \,for all \,$s \in \N$.
As obtained in the proof of Theorem \ref{Sm-Spae dense lineable}, we get for every \,$(n,q) \in \N \times (0,+\infty )$ \,that
$$
\int_\Omega |f_n|^q \, d\mu = \mu (\Omega_{j_1 \cdots j_m}) = {\mu (A) \over 2^{m(n)}}.
$$
Let us denote \,$\gamma := \sum_{{\bf j} \in J} |a_{\bf j}|$. It follows that
\begin{equation*}
\begin{split}
0 &\le \limsup_{n \to \infty} \int_\Omega |Q_n|^q \, d\mu = \limsup_{n \to \infty} \left| \sum_{{\bf j} \in J} a_{\bf j} \cdot e^{\sigma ({\bf j}) \sqrt{m(n)}} \right|
\int_\Omega |f_n|^q \, d\mu \\
&\le \lim_{n \to \infty} \gamma \cdot e^{\sigma ({\bf k}) \sqrt{m(n)}} \cdot {\mu (A) \over 2^{m(n)}} = \gamma \cdot \mu (A) \cdot \lim_{n \to \infty}
  e^{\sigma ({\bf k}) \sqrt{m(n)} - m(n) \log 2} = 0.
\end{split}
\end{equation*}
Consequently, $\lim_{n \to \infty} \int_\Omega |Q_n|^q \, d\mu = 0$, that is, ${\bf Q} \in \bigcap_{q > 0} \mathcal{S}_{L_q}$.

\vskip 3pt

Finally, in the proof of Theorem \ref{Sm-Spae dense lineable} it was shown that, for any given \,$x_0 \in A$, the sequence
\,$(f_n(x_0))$ \,contains infinitely many $0$'s and infinitely many $1$'s. In particular, for infinitely many \,$n \in \N$
\,it holds that
$$
|Q_n(x_0)| = \left| \sum_{{\bf j} \in J} a_{\bf j} \cdot e^{\sigma ({\bf j}) \sqrt{m(n)}} \right| \ge
e^{\sigma ({\bf k}) \sqrt{m(n)}} \cdot \left( |a_{\bf k}| - \sum_{{\bf j} \in J \setminus \{ {\bf k} \}} |a_{\bf j}| \cdot e^{(\sigma ({\bf j}) - \sigma ({\bf k})) \sqrt{m(n)}} \right) ,
$$
and the right hand side tends to $\infty$ because the exponents $\sigma ({\bf j}) - \sigma ({\bf k})$ in the last sum are all negative.
Thus, the sequence $(Q_n(x_0))$ cannot tend to zero, which entails \,${\bf Q} \not \in \mathcal{S}_p$, as desired.
\end{proof}

\begin{theorem} \label{Sm - (Sp u U SLq) algebrable}
Let \,$(\Omega ,\mathcal{A},\mu )$ \,be a nonatomic semifinite measure space. Then the set
\,$\mathcal{S}_m \setminus \Big( \mathcal{S}_p \cup \bigcup_{q > 0} \mathcal{S}_{L_q} \Big)$  \,is strongly $\mathfrak{c}$-algebrable.
\end{theorem}

\begin{proof}
Here and in the next two theorems we shall follow an approach similar to the proof of Theorem \ref{Sm-Spae algebrable}.
Similarly to the construction given in the proof of Theorem \ref{Sm-(Sp u USLq) dense lineable y spaceable}, whose notation we maintain, let us define for each
\,$c > 0$ \,the sequence \,${\bf f}_c = (f_{c,n})$ \,by
$$
f_{c,n} = e^{m(n)^2 c} \cdot f_n,
$$
where \,$(f_n)$ \,is the sequence built
in the proof of Theorem \ref{Sm-Spae dense lineable}, that is,
\,$f_1 = \chi_A$ \,(recall that $0 < \mu (A) < \infty$) \,and \,$f_n = \chi_{\Omega_{j_1 \cdots j_m}}$ \,if \,$n \ge 2$.
Choose a $\Q$-linearly independent set \,$H \subset (0,+\infty )$ \,with ${\rm card} \, (H) = \mathfrak{c}$.
Fix \,$N \in \N$ \,as well as different elements
\,$c_1,c_2, \dots ,c_N$ \,of \,$H$. Also, fix a nonzero polynomial \,$P$ \,in \,$N$ \,variables without constant term, with the form
given in \eqref{Equation Polynomial}, and with \,$J$ \,and the coefficients \,$a_{\bf j}$ \,satisfying the same properties as in the proof of Theorem \ref{Sm-Spae algebrable}.
Recall that the numbers \,$\sigma ({\bf j}) := c_1j_1 + \cdots + c_Nj_N$ (${\bf j} \in J$) \,were pairwise different, and so there is only
one of them, denoted \,$\sigma ({\bf k})$, that is their maximum.

\vskip 3pt

We have to prove that the sequence
$$
{\bf Q} = (Q_n) := P({\bf f}_{c_1}, \dots ,{\bf f}_{c_N})
$$
belongs to \,$\mathcal{S}_m \setminus \Big( \mathcal{S}_p \cup \bigcup_{q > 0} \mathcal{S}_{L_q} \Big)$ (that ${\bf Q} \ne 0$ is automatically implied by this).
Observe first that \,$Q_n = \big( \sum_{{\bf j} \in J} a_{\bf j} \cdot e^{\sigma ({\bf j}) m(n)^2} \big) \cdot f_n$ \,for all \,$n \in \N$.
Now, let us fix \,$(\ve , x_0, q) \in (0,+\infty ) \times A \times (0,+\infty )$. Then we have:
\begin{itemize}
\item Since \,$f_n$ \,is not zero only on \,$\Omega_{j_1 \cdots j_m}$, we get for \,$n \ge 2$ \,that
$$
\phantom{aaaaaa} \mu (|Q_n| > \ve ) \le \mu (Q_n \ne 0) \le \mu (f_n \ne 0) =
\mu (\Omega_{j_1 \cdots j_m}) = 2^{-m(n)} \longrightarrow 0 \hbox{ \ as } \, n \to \infty,
$$
because \,$m(n) \to \infty$. It follows that \,$\lim_{n \to \infty} \mu (|Q_n - 0| > \ve ) = 0$.
Therefore, ${\bf Q} \in \mathcal{S}_m$.
\item Recall, on the one hand, that the sequence \,$(f_n(x_0))$ \,contains infinitely many $1$'s.
On the other hand, note that for such $n$'s we have
\begin{equation*}
\phantom{aaaaaa} |Q_n(x_0)| = \left| \sum_{{\bf j} \in J} a_{\bf j} \cdot e^{\sigma ({\bf j}) m(n)^2} \right|
\ge e^{\sigma ({\bf k}) m(n)^2} \cdot \left( |a_{\bf k}| - \sum_{{\bf j} \in J \setminus \{{\bf k}\} } |a_{\bf j}| \cdot e^{(\sigma ({\bf j}) - \sigma ({\bf k}) ) m(n)^2}                 \right).
\end{equation*}
Since \,$\sigma ({\bf j}) - \sigma ({\bf k}) < 0$ \,for all \,${\bf j} \ne {\bf k}$, we derive that the term \,$\gamma_n$ \,between brackets tends to \,$|a_{\bf k}| > 0$ \,as \,$n \to \infty$
\,and, therefore, $|Q_n(x_0)| \ge \gamma_n \ge |a_{\bf k}|/2$ \,for infinitely many \,$n \in \N$, which prevents \,$(Q_n(x_0))$ \,to tend to zero.
Hence, ${\bf Q} \not \in \mathcal{S}_p$.
\item  According to the previous item, we obtain that \,$|Q_n(x)| \ge e^{\sigma ({\bf k}) m(n)^2} \cdot |a_{\bf k}|/2$ \,for infinitely many
\,$n \in \N$ \,and all \,$x \in A$. Then for such \,$n$'s \,we have
\begin{equation*}
\begin{split}
\int_\Omega |Q_n|^q \, d\mu &= \int_{\Omega_{j_1(n) \cdots j_{m(n)}(n)}} |Q_n|^q \,d\mu \ge e^{\sigma ({\bf k}) m(n)^2 q} \cdot {|a_{\bf k}|^q \over 2^q} \cdot \mu (\Omega_{j_1(n) \cdots j_{m(n)}(n)}) \\
&= {|a_{\bf k}|^q \over 2^q} \cdot 2^{ \sigma ({\bf k}) m(n)^2 q \log 2 - m(n) } \longrightarrow \infty \hbox{ \ as } \, n \to \infty .
\end{split}
\end{equation*}
It follows that \,$(Q_n)$ \,cannot tend to \,$0$ \,in \,$q$-mean.
Thus, ${\bf Q} \not \in \mathcal{S}_{L_q}$.
\end{itemize}
Consequently, our claim about the membership of \,${\bf Q}$ \,holds. The proof is finished.
\end{proof}

\begin{theorem} \label{Spae-Sm algebrable}
Let \,$(\Omega ,\mathcal{A},\mu )$ \,be a measure space satisfying {\rm (P)}. 
Then \,$\mathcal{S}_{p} \setminus \mathcal{S}_m$ \,is strongly $\mathfrak{c}$-algebrable.
In particular, the sets \,
$\mathcal{S}_p \setminus \bigcup_{q > 0} \mathcal{S}_{L_q}$,
$\mathcal{S}_p \setminus \mathcal{S}_{c}$,
$\mathcal{S}_{p} \setminus \mathcal{S}_{au}$ \,and \,$\mathcal{S}_{p} \setminus \mathcal{S}_u$ \,are strongly $\mathfrak{c}$-algebrable.
\end{theorem}

\begin{proof}
By using the construction given in the proof of Theorem \ref{Spae-Sm dense lineable} and thanks to (P), we can start
with a sequence \,$(f_n) \in L_0^{\N}$ \,defined as
$$
f_n = \chi_{A_n},
$$
where \,$(A_n)$ \,is the sequence of measurable sets
satisfying, for some fixed \,$\alpha > 0$, that
$$
\mu (A_n) \ge \alpha \hbox{ \ for all } \, n \in \N .
$$
For each \,$c > 0$, let the sequence \,${\bf f}_c = (f_{c,n})$ \,be given by
$$
f_{c,n} = e^{cn} \cdot f_n.
$$
Again, choose a $\Q$-linearly independent set \,$H \subset (0,+\infty )$ \,with ${\rm card} \, (H) = \mathfrak{c}$, and
consider the family \,$\{ {\bf f}_c = (f_{c,n}) \in L_0^{\N} : \, c \in H \}$. Fix \,$N \in \N$ \,as well as different elements
\,$c_1,c_2, \dots ,c_N$ \,of \,$H$. Also, fix a nonzero polynomial \,$P$ \,in \,$N$ \,variables without constant term, as in \eqref{Equation Polynomial}, and with \,$J$ \,and the coefficients \,$a_{\bf j}$ \,satisfying the same properties as in the proof of Theorem \ref{Sm-Spae algebrable}.
Again, the numbers \,$\sigma ({\bf j})$ \,are pairwise different, and so there is only
one of them, denoted \,$\sigma ({\bf k})$, such that \,$\sigma ({\bf k}) = \max \{ \sigma ({\bf j}) : \, {\bf j} \in J \}$.
Our goal is to show that the sequence
$$
{\bf Q} = (Q_n) := P({\bf f}_{c_1}, \dots ,{\bf f}_{c_N})
$$
belongs to \,$\mathcal{S}_{p} \setminus \mathcal{S}_{m}$.

\vskip 3pt

On the one hand, if we fix \,$x_0 \in \Omega$ \,then, since the \,$A_n$'s \,are pairwise disjoint, $x_0$ belongs to at most one \,$A_n$.
Therefore, there is \,$n_0 \in \N$ \,such that
\,$f_{c_\nu,n}(x_0) = e^{c_\nu n} \cdot \chi_{A_n}(x_0) = 0$ \,for all \,$n > n_0$ \,and all \,$\nu \in \{1, \dots, N \}$.
But \,$P(0,0, \dots ,0) = 0$. Therefore, for all \,$n > n_0$ \,we get
$$
Q_n(x_0) = P(f_{c_1,n}(x_0), \dots ,f_{c_N,n}(x_0)) = 0 \longrightarrow 0 \, \hbox{ as } \, n \to \infty ,
$$
and this holds for all \,$x_0 \in \Omega$. Thus, ${\bf Q} \in \mathcal{S}_{p}$.

\vskip 3pt

On the other hand, observe that for all \,$x \in A_n$ \,we have
$$
|Q_n(x)| = \left| \sum_{{\bf j} \in J} a_{\bf j} \cdot e^{\sigma ({\bf j}) n} \right| \ge
e^{\sigma ({\bf k}) n} \cdot \left( |a_{\bf k}| - \sum_{{\bf j} \in J \setminus \{ {\bf k} \}} |a_{\bf j}| \cdot e^{(\sigma ({\bf j}) - \sigma ({\bf k})) n} \right).
$$
Since the right hand side of the last inequality tends to \,$\infty$ \,as \,$n \to \infty$, we derive the existence of \,$n_0 \in \N$ \,satisfying that
\,$|Q_n(x)| > 1$ \,for all \,$x \in A_n$ \,provided that \,$n > n_0$. This tells us that, if \,$\varepsilon := 1$, then
\,$\mu (|Q_n| > \varepsilon ) \ge \mu (A_n) \ge \alpha$ \,for \,$n$ \,large enough,
and so \,$\mu (|Q_n - 0| > \ve) \not\longrightarrow 0$ \,as \,$n \to \infty$. Consequently, ${\bf Q} \not\in \mathcal{S}_m$, as required.
\end{proof}

\begin{theorem} \label{Sau-Su algebrable}
Let \,$(\Omega ,\mathcal{A},\mu )$ \,be a measure space satisfying {\rm (R)}. 
Then \,$(\mathcal{S}_c \cap \mathcal{S}_{au}) \setminus \mathcal{S}_u$ \,is strongly $\mathfrak{c}$-algebrable.
In particular, $\mathcal{S}_c \setminus \mathcal{S}_u$ \,and
\,$\mathcal{S}_{au} \setminus \mathcal{S}_u$ \,are strongly $\mathfrak{c}$-algebrable.
\end{theorem}

\begin{proof}
We have only to prove the strong $\mathfrak{c}$-algebrability of \,$(\mathcal{S}_c \cap \mathcal{S}_{au}) \setminus \mathcal{S}_u$.
Let us make exactly the same construction (and we keep the notation) as the one in the proof of Theorem \ref{Spae-Sm algebrable}, with the sole exception that the sets
\,$A_n$ \,are replaced \,by measurable sets \,$B_1 \supset B_2 \supset B_3 \supset \cdots$ \,satisfying \,$\mu (B_n) > 0$ \,for all \,$n \in \N$
\,and \,$\lim_{n \to \infty} \mu (B_n) = 0$. This is possible thanks to the hypothesis (R): see the proof of Theorem \ref{Sau-Su dense lineable}.
By extracting a subsequence if necessary, we can assume that \,$\mu (B_n) < 2^{-n}$ \,for all \,$n \in \N$.

\vskip 3pt

Then for \,${\bf Q} = (Q_n) := P({\bf f}_{c_1}, \dots ,{\bf f}_{c_N})$ \,we should prove that \,${\bf Q} \in (\mathcal{S}_{au} \cap \mathcal{S}_{c}) \setminus \mathcal{S}_u$,
where this time
$$
{\bf f}_c = (f_{c,n}) := (e^{cn} \cdot \chi_{B_n}).
$$
For a polynomial \,$P$ \,as in \ref{Equation Polynomial}, we have that
\,$Q_n = \big( \sum_{{\bf j} \in J} a_{\bf j} \cdot e^{\sigma ({\bf j}) n}  \big) \cdot \chi_{B_n}$.
Fix \,$\varepsilon > 0$.

\vskip 3pt

On the one hand, we can select \,$m \in \N$ \,with \,$\mu (B_m) < \varepsilon$.
Let \,$Z_\varepsilon := B_m \in \mathcal{A}$. Then \,$\mu (Z_\varepsilon) < \ve$.
For every \,$x \in \Omega \setminus Z_\varepsilon$ \,we have that \,$x \not\in B_n$ \,for all \,$n \ge m$, and so \,$\chi_{B_n}(x) = 0 = Q_n(x)$ \,for all \,$n \ge m$.
Then
$$
\lim_{n \to \infty} \sup_{x \in \Omega \setminus Z_\varepsilon} |Q_n(x) - 0| = 0,
$$
that is, ${\bf Q} \in \mathcal{S}_{au}$.

\vskip 3pt

On the other hand, notice that since \,$\{|Q_n| > \ve \} \subset \{Q_n \ne 0 \} \subset B_{n}$, we get
$$
\sum_{n=1}^{\infty} \mu (\{|Q_n - 0| > \ve \}) \le \sum_{n=1}^{\infty} \mu (B_n) \le \sum_{n=1}^{\infty} 2^{-n} = 1 < \infty .
$$
From this we infer that ${\bf F} \in \mathcal{S}_{c}$.

\vskip 3pt

Finally, fix \,$Z \in \mathcal{A}$ \,with \,$\mu (Z) = 0$. Since \,$\mu (B_m) > 0$, the intersection \,$(\Omega \setminus Z ) \cap B_n$
\,cannot be empty. It follows the existence of a sequence \,$(x_n) \subset \Omega$ \,with \,$x_n \in \Omega \setminus Z$ \,and \,$\chi_{B_n} (x_n) = 1$
\,for all \,$n \in \N$. Consequently,
$$
\sup_{x \in \Omega \setminus Z} |Q_n(x)-0| \ge |Q_n(x_n)| = \left| \sum_{{\bf j} \in J} a_{\bf j} \cdot e^{\sigma ({\bf j}) n} \right|
\ge e^{\sigma ({\bf k}) n} \cdot \left( |a_{\bf k}| - \sum_{{\bf j} \in J \setminus \{{\bf k}\} } |a_{\bf j}| \cdot e^{(\sigma ({\bf j}) - \sigma ({\bf k}) ) n}                  \right).
$$
Since the right hand side of the last chain of inequalities tends to $\infty$ as $n \to \infty$, we derive
that \,$(Q_n)$ \,cannot tend uniformly a.e., and so \,${\bf Q} \not\in \mathcal{S}_u$. The proof is finished.
\end{proof}

\begin{theorem} \label{Su-SL algebrable}
Let \,$(\Omega ,\mathcal{A},\mu )$ \,be a measure space. The following properties are equi\-va\-lent:
\begin{enumerate}
\item[\rm (a)] $\mu (\Omega ) = \infty$.
\item[\rm (b)] $\mathcal{S}_u \setminus \bigcup_{q > 0} \mathcal{S}_{L_q} \ne \varnothing$.
\item[\rm (c)] $\mathcal{S}_u \setminus \bigcup_{q > 0} \mathcal{S}_{L_q}$ \,is strongly $\mathfrak{c}$-algebrable.
\end{enumerate}
In particular, if \,$\mu (\Omega ) = \infty$, then the sets
\,$\mathcal{S}_p \setminus \bigcup_{q > 0} \mathcal{S}_{L_q}$,
$\mathcal{S}_m \setminus \bigcup_{q > 0} \mathcal{S}_{L_q}$ \,and \,
$\mathcal{S}_{au} \setminus \bigcup_{q > 0} \mathcal{S}_{L_q}$ \,are strongly $\mathfrak{c}$-algebrable.
\end{theorem}

\begin{proof}
Recall that, within a finite measure space, uniform convergence a.e.~implies convergence in $q$-mean. Then, trivially, we have (c) $\Longrightarrow$ (b) $\Longrightarrow$ (a).
It remains to prove (a) $\Longrightarrow$ (c).

\vskip 3pt

To this end, select a $\Q$-linearly independent set \,$H \subset (0,+\infty )$ \,with \,${\rm card} \, (H) = \mathfrak{c}$, and consider the family
\,$\{ {\bf f}_c : \, c \in H \}$ \,of constant functions given by
$$
{\bf f}_c = (e^{-cn}).
$$
Fix a nonzero polynomial \,$P$ \,without constant term as the one given by \eqref{Equation Polynomial} in the proof of Theorem \ref{Sm-Spae algebrable}. We keep all notations and conditions concerning to \,$P$ \,in the mentioned proof, with the sole exception that this time
\,${\bf k}$ \,represents the unique $N$-index in \,$J$ \,such that
\,$\sigma ({\bf k}) = \min \{ \sigma ({\bf j}) : \, {\bf j} \in J \}$. Then
 \,$\sigma ({\bf k}) - \sigma ({\bf j}) < 0$ \,for all \,${\bf j} \in J \setminus \{ {\bf k} \}$.
Let us define
$$
{\bf Q} = (Q_n) := P({\bf f}_{c_1}, \dots , {\bf f}_{c_N}),
$$
that is, ${\bf Q}$ is the sequence whose terms are the constant functions
$$
Q_n = \sum_{{\bf j} \in J} a_{\bf j} \cdot e^{- \sigma ({\bf j}) n} \quad (n=1,2,\dots ).
$$
It should be shown that \,${\bf Q} \in (\mathcal{S}_u \setminus \bigcup_{q > 0} \mathcal{S}_{L_q}) \setminus \{0\} = \mathcal{S}_u \setminus \bigcup_{q > 0} \mathcal{S}_{L_q}$.
Firstly, observe that
$$
0 \le \sup_{x \in \Omega} |Q_n(x) - 0| = \left|  \sum_{{\bf j} \in J} a_{\bf j} \cdot e^{- \sigma ({\bf j}) n} \right|
\le \left(  \sum_{{\bf j} \in J} |a_{\bf j}| \right) \cdot e^{- \sigma ({\bf k}) n} \longrightarrow 0 \hbox{ \ as } \,n \to \infty ,
$$
which tells us that \,${\bf Q} \in \mathcal{S}_u$.

\vskip 3pt

Finally, notice that since \,$\mu (\Omega ) = \infty$ \,and each \,$Q_n$ \,is constant, then we shall have for all \,$q > 0$ \,that \,$\int_\Omega |Q_n|^q \,d\mu = \infty$ \,as soon as \,$Q_n \ne 0$. If this were asymptotically true, then one would obtain \,$\int_\Omega |Q_n|^q \,d\mu \to \infty$ \,as \,$n \to \infty$, that would imply
\,${\bf Q} \not \in \bigcup_{q > 0} \mathcal{S}_{L_q}$, as required. Since
$$
|Q_n| \ge e^{- \sigma ({\bf k}) n} \cdot \bigg( |a_{\bf k}| - \sum_{{\bf j} \in J \setminus \{ {\bf k} \} } |a_{\bf j}| \cdot e^{(\sigma ({\bf k}) - \sigma ({\bf j})) n} \bigg)
$$
and the sum inside the right hand side tends to \,$0$ \,as \,$n \to \infty$, we get \,$|Q_n| \ge e^{- \sigma ({\bf k}) n} \cdot {|a_{\bf k}| \over 2} > 0$ \,for \,$n$ \,large enough. This concludes the proof.
\end{proof}

\section{Final remarks and questions}

\noindent {\bf 1.} Many of the results contained in (A), (B), (C), (D), (E) and (F) of section 2 have been highly improved/extended/covered by our Theorems
\ref{Sm-Spae dense lineable} to \ref{Su-SL algebrable}. This holds because: $([0,1], \mathcal{B},\lambda )$ \,is a nonatomic probability space satisfying (S) and (R);
$([0,+\infty] ), \mathcal{B},\lambda )$ \,is a nonatomic $\sigma$-finite measure space satisfying (S), (P), (Q), (R) and $\lambda ([0,+\infty ) ) = \infty$; every probability space is $\sigma$-finite; every $\sigma$-finite measure space is semifinite; every nonatomic semifinite measure space satisfies (R); and every nonatomic semifinite measure space with infinite measure satisfies (P) and (Q).

\vskip 4pt

\noindent {\bf 2.} As said in sections 1 and 2, our goal in this paper has been to deal with some kinds of convergence {\it to the zero function.} Nevertheless, it is fair to say that other interesting convergence --not necessarily to zero-- properties under the point of view of li\-ne\-a\-bi\-li\-ty appear in the papers \cite{vecina,polacos2023,conejerofms,fesetu}. For instance: convergence in distribution and in variation, together with convergence of sums of random variables on a probability space and convex-lineability properties, are considered in \cite{vecina};
also in the probability theory setting, Conejero {\it et al.}~\cite{conejerofms} studied some lineability and algebrability pro\-blems, as for example, convergent not $L_1$-unbounded martingales, pointwise convergent random variables whose means do not converge to the expected value, and stochastic processes that are $L_2$-bounded and convergent but not pointwise convergent in a null set; still in a (nonatomic) probability space, lineability and algebrability are discovered
in \cite{polacos2023} (see also \cite{polacos2018}) for several fa\-mi\-lies of \,$L_0^{\N}$ \,of independent random variables with additional properties connected with various types of convergence, laws of large numbers, non-preservation of convergence in measure under C\`esaro means, and Markov and Kolgomorov conditions; in \cite{fesetu}, it is proved that when omitting one condition in classical convergence results (such as the dominated convergence theorem, Fatou's Lemma, or the Strong Law of Large Numbers), it is possible to find lineability, algebrability, coneability or latticeability for the fa\-mi\-lies of the respective counterexamples and, in addition, it is shown that in \,$([0,1], \mathcal{B},\lambda )$ \,the fa\-mi\-lies of sequences of random variables converging in probability but (i) not converging outside a set of measure zero or (ii) not converging in arithmetic mean are also linearly very large. Moreover, sequences of continuous (on a locally compact Hausdorff space) $L_1$-undominated functions tending to zero in $1$-mean and uniformly are studied in \cite{bercalmur2020} (see also \cite{bernalordonez2014}).

\vskip 4pt

\noindent {\bf 3.} In the present paper, we have put the emphasis on six modes of convergence (pointwise a.e, in measure, in $q$-mean, uniform a.e, almost uniform, complete).
It would be interesting to complement the study by introducing other kinds of convergence to be compared with the previous ones or between them, such as
convergence in distribution and convergence in variation, which are more natural in the setting of probability spaces.

\vskip 4pt

\noindent {\bf 4.} In the specific framework of Theorems \ref{Sm-Spae dense lineable} and \ref{Su-SL dense lineable}, and under appropriate assumptions on the measure space,
the following questions arise naturally: \hfil\break
(1) Is \,$\displaystyle{ \Big( \bigcap_{q > 0} \mathcal{S}_{L_q} \Big)} \setminus \mathcal{S}_{p}$ \,spaceable in \,$L_0^{\N}$?  \hfil\break
(2) Is \,$\mathcal{S}_u \setminus \bigcup_{q > 0} \mathcal{S}_{L_q}$ \,spaceable in \,$L_0^{\N}$?


%
%
%
%
%
%
%


\end{document}